\documentclass[10pt,letterpaper]{article}
\usepackage[top=1in,left=1.2in,right=1.2in]{geometry}

\usepackage{changepage}

\usepackage[utf8]{inputenc}

\usepackage{textcomp,marvosym}

\usepackage{fixltx2e}

\usepackage{amsmath,amssymb}

\usepackage{cite}

\usepackage{nameref,hyperref}

\usepackage[right]{lineno}

\usepackage{microtype}
\DisableLigatures[f]{encoding = *, family = * }

\usepackage{rotating}

\raggedright
\usepackage[aboveskip=1pt,labelfont=bf,labelsep=period,justification=raggedright,singlelinecheck=off,font=small]{caption}

\makeatletter
\renewcommand{\@biblabel}[1]{\quad#1.}
\makeatother

\date{}

\usepackage{lastpage,fancyhdr,graphicx}
\usepackage{epstopdf}
\usepackage{etex} % allows using more packages - otherwise compilation errors

\usepackage{url}
\usepackage{algorithm}
\usepackage{algorithmic}
\usepackage{amsmath}
\usepackage{amsfonts}
\usepackage{amssymb}
\usepackage{amsthm}
\usepackage{colortbl}
\usepackage{floatflt} % for floating inlets
\usepackage{graphicx}
\usepackage{grffile} % allows extra "." in file names
\usepackage{latexsym}
\usepackage{multicol} % for multiple column lists
\usepackage{multirow}
\usepackage{paralist}
\usepackage{pgfplotstable}
\usepackage{relsize} % for \mathlarger
\usepackage{subcaption}
\usepackage[titletoc]{appendix}
\usepackage[textsize=tiny,english,colorinlistoftodos,textwidth=5cm]{todonotes}
\usepackage{varwidth}
\usepackage{wrapfig}
\usepackage{pdflscape} % for rotating single pages
\usepackage{pgfplotstable} % for CSV based table
\usepackage[para]{footmisc} % for footnotes in one line

\usepackage{tikz}
\usetikzlibrary{plotmarks}
\usepackage{longtable}

\usepackage{ctable}
\usepackage[normalem]{ulem}

\newcommand{\R}{{\mathbb R}}
\DeclareMathOperator*{\argmax}{arg\,max}
\newcommand\abs[1]{\left|#1\right|}

\theoremstyle{definition}
\newtheorem{algo}{Algorithm}
\newcommand{\algobox}[3][]{\vspace{.5\baselineskip}
\fbox{\parbox{#2}{\vspace{-.5\baselineskip}
\begin{algo}[#1]\leavevmode
#3
\end{algo}\vspace{-.5\baselineskip}}}
\vspace{\baselineskip}
}

\newtheorem*{remark*}{Remark}

\DeclareMathOperator{\sign}{sign}
\DeclareMathOperator{\supp}{supp}
\colorlet{LightRed}{white!70!red}
\colorlet{LightBlue}{white!70!blue}
\colorlet{LightGreen}{white!70!green}
\definecolor{DarkGreen}{HTML}{009600}

\newcommand\reduline{\bgroup\markoverwith{\textcolor{red}{\rule[0.5ex]{2pt}{0.4pt}}}\ULon}

\begin{document}
% \vspace*{0.35in}

% Please capitalize all terms in the title except conjunctions, prepositions, and articles.
\begin{flushleft}
{\Large
\textbf\newline{Sparse Proteomics Analysis -- A compressed sensing-based approach for feature selection and classification of high-dimensional proteomics mass spectrometry data}
}
\newline
\\
% change footnote symbols to
\renewcommand*{\thefootnote}{\alph{footnote}}
Tim OF Conrad\textsuperscript{1,*},
Martin Genzel\textsuperscript{2}\footnote{genzel@math.tu-berlin.de},
Nada Cvetkovic\textsuperscript{1}\footnote{ncvetkovic@math.fu-berlin.de},
Niklas Wulkow\textsuperscript{1}\footnote{nwulkow@math.fu-berlin.de},
Alexander Leichtle\textsuperscript{3}\footnote{alexander.leichtle@insel.ch},
Jan Vybiral\textsuperscript{4}\footnote{vybiral@karlin.mff.cuni.cz},
Gitta Kutyniok\textsuperscript{2}\footnote{kutyniok@math.tu-berlin.de},
Christof Sch\"utte\textsuperscript{1,5}\footnote{schuette@math.fu-berlin.de},
% change back footnote symbold
\renewcommand*{\thefootnote}{\arabic{footnote}}
\setcounter{footnote}{0}
\\
\bigskip
\bf{1} Department of Mathematics, Freie Universit\"at Berlin, Germany
\\
\bf{2} Department of Mathematics, Technische Universit\"at Berlin, Germany
\\
\bf{3} Center of Laboratory Medicine, Inselspital - Bern University Hospital, Switzerland
\\
\bf{4} Department of Mathematical Analysis, Charles University, Prague, Czech Republic
\\
\bf{5} Zuse Institute Berlin, Berlin, Germany
\\
\bigskip

% Use the asterisk to denote corresponding authorship and provide email address in note below.
* Corresponding author: conrad@math.fu-berlin.de; Postal address: Department of Mathematics, Freie Universit\"at Berlin, Arnimallee 6, 14195 Berlin

\end{flushleft}

\section*{Abstract}
\paragraph{Background}
High-throughput proteomics techniques, such as mass spectrometry (MS)-based approaches, produce very high-dimensional data-sets. In a clinical setting one is often interested in how mass spectra differ between patients of different classes, for example spectra from healthy patients vs. spectra from patients having a particular disease. Machine learning algorithms are needed to (a) identify these discriminating features and (b) classify unknown spectra based on this feature set. Since the acquired data is usually noisy, the algorithms should be robust against noise and outliers, while the identified feature set should be as small as possible.

\paragraph{Results}
We present a new algorithm, \emph{Sparse Proteomics Analysis} (\emph{SPA}), based on the theory of compressed sensing that allows us to identify a minimal discriminating set of features from mass spectrometry data-sets. We show (1) how our method performs on artificial and real-world data-sets, (2) that its performance is competitive with standard (and widely used) algorithms for analyzing proteomics data, and (3) that it is robust against random and systematic noise. We further demonstrate the applicability of our algorithm to two previously published clinical data-sets.

%%%%%%%%%%%%%%%%%%%%%%%%%%%%%%%%%%%%%%%%%%%%%%%%%%%%%%%%%%%%%%%%%%%%%%%%%%%%%%%%
%%%%%%%%%%%%%%%%%%%%%%%%%%%%%%%%%%%%%%%%%%%%%%%%%%%%%%%%%%%%%%%%%%%%%%%%%%%%%%%%
\section*{Availability of data and methods}
The method source-code can be downloaded from our homepage: \url{http://software.medicalbioinformatics.de}. The used data will be made available through BioMed Central.

%%%%%%%%%%%%%%%%%%%%%%%%%%%%%%%%%%%%%%%%%%%%%%%%%%%%%%%%%%%%%%%%%%%%%%%%%%%%%%%%
%%%%%%%%%%%%%%%%%%%%%%%%%%%%%%%%%%%%%%%%%%%%%%%%%%%%%%%%%%%%%%%%%%%%%%%%%%%%%%%%
\section*{List of abbreviations}
AIC - Akaike Information Criterion; BIC - Bayesian Information Criterion; CS - Compressed Sensing; MALDI-TOF: Matrix-Assisted Laser Desorption Ionization Time-Of-Flight; ML - Machine Learning; MS - Mass Spectrometry;  SPA - Sparse Proteomics Analysis; SVM - Support Vector Machine; TP - True Positive; TN - True Negative; FP - False Positive; FN - False Negative

% \linenumbers

%%%%%%%%%%%%%%%%%%%%%%%%%%%%%%%%%%%%%%%%%%%%%%%%%%%%%%%%%%%%%%%%%%%%%%%%%%%%%%%%
%%%%%%%%%%%%%%%%%%%%%%%%%%%%%%%%%%%%%%%%%%%%%%%%%%%%%%%%%%%%%%%%%%%%%%%%%%%%%%%%

%%%%%%%%%%%%%%%%%%%%%%%%%%%%%%%%%%%%%%%%%%%%%%%%%%%%%%%%%%%%%%%%%%%%%%%%%%%%%%
\section{Introduction}
%%%%%%%%%%%%%%%%%%%%%%%%%%%%%%%%%%%%%%%%%%%%%%%%%%%%%%%%%%%%%%%%%%%%%%%%%%%%%%

During the last decade, high-throughput assays systems\footnote{Assays, e.g. immunoassays, are used in molecular diagnostics to detect concentrations of specific molecules even in low concentrations from a biological sample, such as blood\cite{Rissin2010}.} for measuring a variety of different biological sources have become standard in modern laboratories. This allows for the quick and cheap creation of very large data-sets which characterize for example the status of a cell by its billions of constituents, e.g. nucleotides, RNAs,  contained proteins, or metabolites. Ideally, analyzing these massive data-sets leads to a better understanding of the underlying biological processes. Especially in the context of characterizing---and ultimately understanding---diseases, a first step is often to find significant differences in the data between samples from healthy and diseased individuals. There are many successful examples where this approach based on -omics data (e.g., genomics, proteomics, or metabolomics) led to the identification of biological markers, enabling a new type of molecular diagnostics. We call a collection of biological markers that represents the differences on the data level a \emph{disease fingerprint}.

Many disease-relevant mechanisms are controlled by proteins (e.g. hormones), which can be detected in biological samples (blood, urine, etc.) using \emph{mass spectrometry} (\emph{MS}). This technique allows (potentially) for monitoring the entire set of proteins---the so-called proteome---in a given sample. Due to its wide availability in hospitals, MS-based proteomics can bring the next wave of progress in diagnostics, since even subtle changes in the proteome can be detected and linked to disease onset and progression \cite{Aebersold2003,Petricoin2006,Rai2004,Coombes2005}.

\begin{figure}

\includegraphics[width=\textwidth]{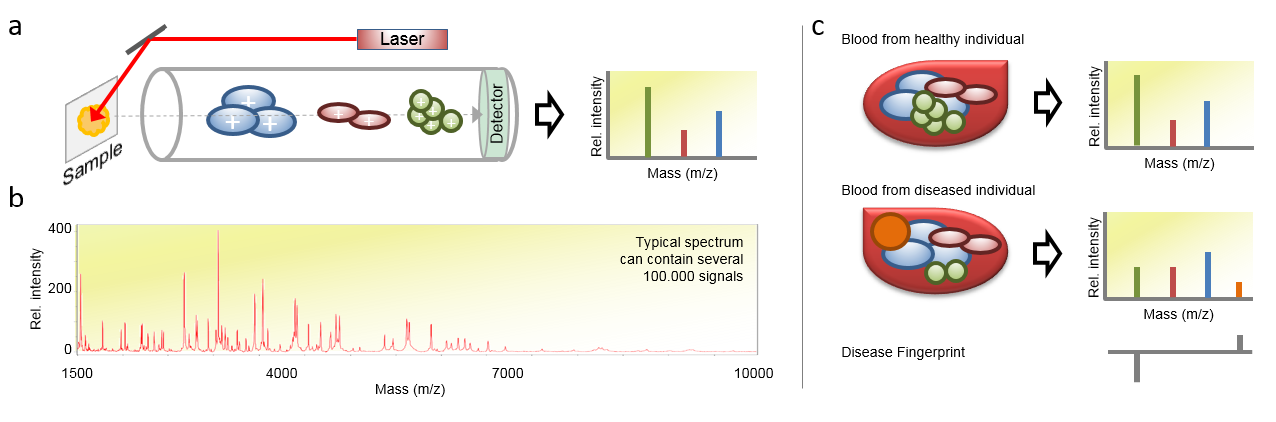}
\caption{(a) Schematic outline of a linear matrix-assisted laser desorption ionization (MALDI)--time-of-flight (TOF) mass spectrometer (MS). During the measurement process, the molecules of the examined sample are ionized, vaporized and finally analyzed by their respective time-of-flight through an electric field. This process generates a plot (mass spectrum) having mass-to-charge ratio ($m/z$) on the $x$-axis and intensity (ion count) on the $y$-axis. (b) Typical mass spectrum for a mass range of 1500--10.000 Dalton. (c) Example of a disease fingerprint, created by comparing mass spectra from a healthy and a diseased individual.}
\label{fig1}
\end{figure}

%%%%%%%%%%%%%%%
\vspace{0.3cm}
\textbf{Disease Fingerprints:}
The main idea of the identification of \emph{disease fingerprints} using MS-based proteomics is sketched in Fig. \ref{fig1}:

(a) A mass spectrum is generated reflecting the constitution of a given (blood-)sample with respect to contained molecules. (b) Based on mass spectra from two sample groups (representing a healthy control group and a group having a particular disease) differences are detected. This set of differences precisely corresponds to a \emph{disease fingerprint}, since it represents a trace caused by a particular disease in the proteome. Several studies have shown that this approach works well in practice and found differences do indeed reflect correlations between changes in the mass spectrum, the proteome, and phenotypic changes (\cite{Liotta2003, Phizicky2003, Issaq2007,Stuehler2004,Sitek2012}). Panels of proteomic markers (fingerprints) have been shown to be more sensitive and specific than conventionally biomarker approaches \cite{Petricoin2006}, for example when diagnosing cancer \cite{Fiedler2009,Strenziok2009,Leichtle2011}. However, a single proteomics data-set can contain tens of millions of signals which is many orders of magnitudes larger than the number of available observations in a typical study.

Our ultimate goal is therefore to build a library of proteomics disease fingerprints which are extracted from high-throughput MS experiments. These would enable to diagnose diseases based on their proteomic fingerprints---just by analyzing an individual's proteome. Ideally, these fingerprints are of low-complexity allowing easy interpretation by experts, e.g. medical doctors, and the implementation of medical assays for routine diagnostics, e.g. in an hospital environment. Clearly, the less components an assay is composed of, the easier it is to implement and interpret.

Thus, a fingerprint should only consist of a minimal collection of proteins specific for a particular disease and should be robust against noisy measurements. On the other hand, the acquired data from the high-throughput experiments is very high-dimensional and contains large amounts of random and systematic noise which makes an automatic analysis of mass spectra a very challenging task. Hence, the discovery of biomarkers is still a widely open research topic and there are several analytic problems that hinder reproduction of results (see \cite{Diao2011} for example).

Despite these challenges there is indeed hope that these disease specific, low-complexity fingerprints exist: It has been shown for several cancer types that a small numbers of genes and proteins can be identified that serve as biomarkers (e.g. for lung cancer \cite{Marrugal2016}, breast cancer \cite{Tang2016} or pancreas cancer \cite{Le2016}). This means that only a few signals in a mass spectrum can be used to derive a sparse classifier.

%%%%%%%%%%%%%%%
\vspace{0.3cm}
\textbf{MS1 Data:}
In this work we consider mass spectrometry (MS) data acquired from a standard MALDI-TOF instrument because it is easy to obtain using comparatively cheap MS-instruments which are widely available, e.g. in hospitals.
Opposed to other approaches such as tandem mass spectrometry (MS/MS), we directly work on the raw data acquired in \emph{profile mode} and do not aim for identification. Thus, each mass spectrum (sample) always has the same number of $d$ dimensions (number of entries). \footnote{The data-sets used in this paper contain $d=42.381$ dimensions in each MS1 spectrum but our approach is not limited by that.} Recall, that the entries in a mass spectrum are a weight-ordered list of ion-counts of the respective ion-masses. (See also Fig. \ref{fig1}.)

One of the reasons for this is that standard approaches for MS data analysis usually convert the MS data to peak lists as a first step and work on the converted data. However, signals can be missed by this conversion step due to noise or missing values in the raw data which hinders peak detection. Opposed to this, our approach does not rely on any peak identification but works on the raw data. This allows for a more robust analysis in presence of noise which is a typical challenge in MS data analysis.

%%%%%%%%%%%%%%%%%%%%%%%%%%%%%%%%%%%%%%%%%%%%%%%%%%%%%%%%%%%%%%%%%%%%%%%%%%%%%%
\subsection{Problem Definition}

%g\vspace{0.1cm}

In this article, we will focus on the following problem setting:

% \vspace{0.2cm}

We assume that we are given data of $n$ mass spectra derived from $n$ biological samples (e.g. from blood of $n$ individual patients) in form of $n$ pairs $\{(x_i,y_i)\}_{i{=}1{\ldots}n}$. Here, $x_i \in \R^d$ represents the mass spectrum of the $i$-th sample (e.g. the $i$-th patient) and $y_i \in \{-1,+1\}$ its respective class, e.g., healthy or diseased. Thus, each $x_i$ (representing an individual mass spectrum) contains $d$ entries. 
The goal is to identify a (small) set of features, i.e. indices in the mass spectrum, separating these two classes. Thus, a feature represents a specific position (or mass) in a mass spectrum in which the two groups (e.g. healthy vs. diseased) differ.
This corresponds to the well known problem of \emph{feature selection}\footnote{In feature selection, one is interested in identifying relevant dimensions of the data (features) which can be used to distinguish between two (or more) classes within a data-set.} and leads to a potential disease fingerprint for the given data.

Mathematically, this can be formulated as the identification of a \emph{feature vector} $\omega_0 = (\omega_{0,1}, \dots, \omega_{0,d}) \in \R^d$ such that\footnote{Here, $\sign(\cdot)$ denotes the sign function, i.e., $\sign(t) = 1$ if $t \geq 0$ and $\sign(t) = -1$ if $t<0$.}
\begin{equation}\label{eq:binarymodel}
y_i = \sign(f_{\omega_0}(x_i)) \quad \text{for ``many'' samples $i=1,\ldots,n$,}
\end{equation}
with a \emph{linear decision function} $f_{\omega_0}(x_i) := \langle {\omega_0}, x_i \rangle = \sum_{j=1}^d \omega_{0,j} x_{i,j}$.

 From a geometric perspective, this means that the hyperplane with normal vector $\omega_0$ appropriately separates the data-points of the respective classes.
% Then the predicted class for a given spectrum $x_i$ depends on the sign of $f_\omega(x_i)$.
This means that $\omega_0$ can be used as a linear classifier where each entry of $\omega_0$ corresponds to a specific position in a spectrum and the non-zero entries (which we call features) indicate their significance. Our goal is therefore to learn a sparse $\omega_0$ for which Eqn. \ref{eq:binarymodel} holds. As a particular consequence, a classifier based on such $\omega_0$ will yield good prediction accuracy. 

In most realistic scenarios for feature selection, unfortunately, the number of features is much larger than available samples ($d \gg n$) and the data suffers from noisy measurements.
For these reasons, the number of feasible classifiers $\omega_0$ can become extremely large, so that the problem of \emph{overfitting} can occur.
In order to allow interpretability and generalization of the classifier, it is in fact inevitable to restrict the solution space for $\omega_0$.
In this paper, we focus on very \emph{sparse}\footnote{We call a vector sparse if the number of non-zero entries is small.} vectors $\omega_0$ satisfying \eqref{eq:binarymodel}, which precisely reflects our wish for a minimal disease fingerprint.

At this point, it should be emphasized that \eqref{eq:binarymodel} does not need to hold for \emph{all} samples but rather for most of them.
Allowing for such a small ``mismatch'' in the model, we incorporate the crucial fact that a simple binary output model, such as \eqref{eq:binarymodel}, might describe the disease label only with high accuracy but not necessarily exactly.
In turn, this asks for a certain robustness of the used method against wrong predictions with regard to \eqref{eq:binarymodel}.

We will approach this challenge by formulating the feature selection problem as a constrained (or regularized) optimization problem:
\begin{equation} \label{eqn:featuresel_intro}
\min_{\omega \in \R^d} \sum_{i=1}^n L( y_i, f_\omega(x_i) ) \quad \text{subject to $R (\omega) \leq \lambda$},
\end{equation}
where $L\colon \R \times \R \to \R$ is a \emph{loss} (error) function, $R\colon \R \to \R$ is a \emph{regularization} (cost) function that encourages a particular structure of $\omega$ (e.g., sparsity), and the parameter $\lambda \ge 0$ controls the degree of model complexity. Given any potential feature vector $\omega$ and the (true) output label $y$, the loss function $L(y,f_{\omega}(x))$ measures the discrepancy between the actual and the desired prediction.

\vspace{0.3cm}
\noindent As already pointed out, we are particularly interested in a method that produces \emph{optimal and robust solutions} in the following situation:

\begin{itemize}
\item  The input data $(x, y)$ are noisy,
\item  the number of data dimensions $d$ is large (typically: $d = 10^5 \ldots 10^8$),
\item  the number of samples $n$ is relatively small (typically: $n = 10^2 \ldots 10^4$), and
\item  the set of highly-relevant features is small (i.e., a minimal disease fingerprint indeed exists), which corresponds to a small number of non-zero elements in $\omega_0$ (typically: $\#\{i \mid \omega_{0,i} \neq 0\} \ll 100$).
\end{itemize}

\vspace{0.3cm}
\noindent On the contrary, we are not mainly interested in the methods' overall classification performance. Measures of classification performance such as accuracy are indicators whether a learned classifier accurately separates the data into classes. In our case, we assume that the data can be characterized well by a \emph{sparse} classifier $\omega_0$ whose non-zero entries are those used for classification and are therefore of medical relevance. That means, if $\omega_0$ is sparse and leads to good classification accuracy then only a few entries contribute and medical interpretation becomes feasible. 
However, if there does not exists a sparse $\omega_0$ such that Eqn. \ref{eq:binarymodel} holds, there is strong evidence that no sparse (simple) characterization is possible. This indicates that the underlying biological mechanisms are too complex to be captured by a sparse (simple) model. If this is the case, every sparsity-encouraging method will fail, meaning that a sparse classifier will always give poor classification. As a consequence, an important assumption of this work is that a sparse $\omega_0$ (ground-truth) exists. 

As we will see later it is often possible to find \emph{non}-sparse classifiers which achieve better classification accuracy. This might be favorable in some situations in which the main focus is indeed on overall classification accuracy. However, in these situations overfitting becomes an issue and the identification of interpretable, highly-discriminative features might be extremely difficult. In the context of MS-data analysis such a classifier would be especially hard to interpret because of the very high dimensionality of the data.

%%%%%%%%%%%%%%%%%%%%%%%%%%%%%%%%%%%%%%%%%%%%%
\subsection{State of the Art in Sparse Feature Selection}
There are numerous approaches for feature selection which mainly fall into three categories:
\begin{itemize}
\item \textbf{Filters}: Using some score or correlation function (e.g., based on Fisher's, t-test, information theoretic criteria) evaluating the importance of each feature in a \emph{univariate} way and taking the top-rated features.
\item \textbf{Wrappers}: Using machine-learning algorithms to evaluate and choose features using some search strategy (e.g. simulated annealing or genetic algorithms).
\item \textbf{Embedded methods}: Selecting variables by directly optimizing an objective function (usually in a multivariate way) with respect to: goodness-of-fit and (optionally) number of features. This could be achieved with algorithms like least-square regression, support vector machines (SVM), or decision trees.
\end{itemize}
In this paper, we will mainly focus on \emph{embedded methods}. Regarding this category, the literature contains several well-known options for choosing combinations of loss and regularization functions (cf. \eqref{eqn:featuresel_intro}), some of which are exemplarily listed in Table~\ref{tab::loss-and-reg-overview}.

\begin{table}
\begin{center}
\vspace{0.5cm}
\begin{tabular}{ l || c | c }
  \hline
  Name & Loss function ($L$) & Regularizer ($R$) \\
  \hline
  AIC/BIC & $\| y - \langle \omega,x \rangle \|_2$ & $\| \omega \|_0$ \\
  Lasso & $\| y - \langle \omega,x \rangle \|_2$ & $\| \omega \|_1$ \\
  Elastic Net & $\| y - \langle \omega,x \rangle \|_2$ & $\| \omega \|^2_2$ + $\| \omega \|_1$ \\
  Regularized Least Absolute  & & \\
  Deviations Regression & $\| y - \langle \omega,x \rangle \|_1$ & $\| \omega \|_1$ \\
  Classic SVM & $\max(0, 1 - y \langle \omega,x \rangle )$\textsuperscript{*} & $ \frac{1}{2} \| \omega \|^2_2$ \\
  $\ell_1$-SVM & $\max(0, 1 - y \langle \omega,x \rangle )$\textsuperscript{*} & $ \frac{1}{2} \| \omega \|_1$ \\
  Logistic Regression & $\log( 1+ \exp(-y \langle \omega,x \rangle ) )$ & $ \frac{1}{2} \| \omega \|_1$ \\
  \hline
  \multicolumn{3}{l}{\textsuperscript{*}\footnotesize{This is the so called \emph{Hinge loss}.}}
\end{tabular}
\end{center}
\caption{Prominent options for choosing loss function and regularizer in feature extraction algorithms. The $\ell_{1}$- and $\ell_2$-norm of a vector $z = (z_1, \ldots, z_d) \in \R^d$ are defined by $\|z\|_1=\sum_{j=1}^d|z_i|$ and $\|z\|_2=(\sum_{j=1}^d |z_i|^2)^{1/2}$, respectively. The ``$\ell_0$-norm'' $\|z\|_0$, simply counts the number of non-zero entries of $z$.}
\label{tab::loss-and-reg-overview}
\end{table}

\begin{figure}[!ht]
% intro_example.png
\includegraphics[width=\textwidth]{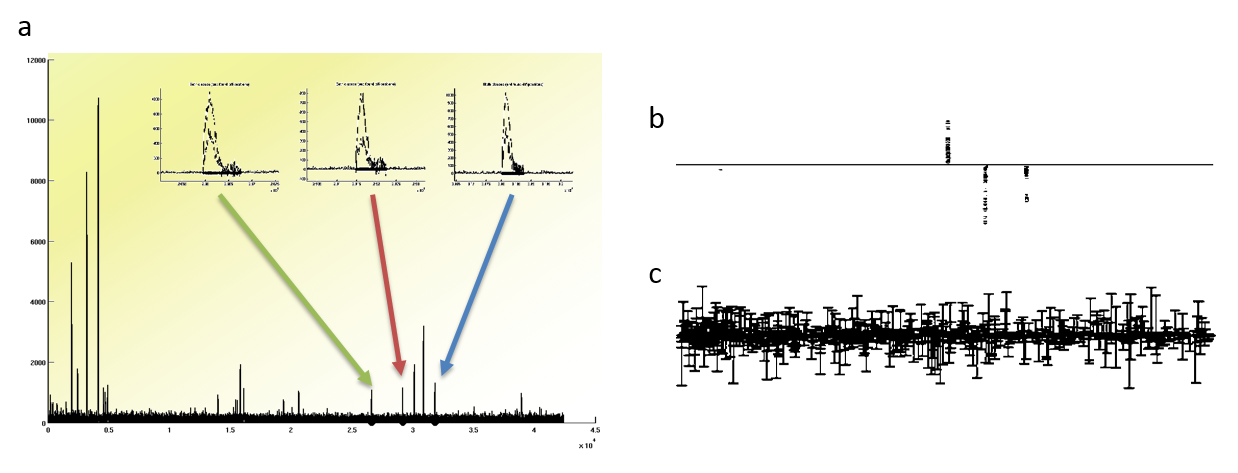}
\caption{(a) Overlaid spectra from two different groups. The three peaks marked by the arrows (magnified in the inlays) represent the underlying differences between the two groups. (b) Sparse $\omega$ found by a $\ell_1$-regularized method ($\ell_1$-SVM). (c) $\omega$ found by $\ell_2$-regularized method (classical SVM).}
% fig:intro_example
\label{fig2}
\end{figure}

Different combinations can influence the results dramatically: Fig.~\ref{fig2} demonstrates the effect of sparsity by comparing a $\ell_2$- and $\ell_1$-regularized version.
In this example, a proteomics data-set was created that contains three discriminant features between the two sub-groups. It can be easily seen how the results differ: While the $\ell_1$-based result is optimized for selecting only a few features, the $\ell_2$-variant selects much more features which in turn results in a better fit of the observation model. In this paper, we are interested in developing a method that selects as few features as possible while achieving the best possible fit under this constraint.
This is in contrast to methods that aim at only achieving the best possible fit. A low-complexity model is of particular interest in biological applications because each selected feature is usually analyzed in subsequent experiments, which creates additional costs.

Various approaches can be used to assess the outcome $\omega$ of a feature selection method, when appropriate training and test data are available. We will use the following three measures of quality: (i) correctness of the selected features, (ii) size of the selected feature set, (iii) performance of classifying an unknown test set (specificity, sensitivity, accuracy). Obviously, (i) can only be used if the correct features are known, which is the case in our benchmark data-sets (for more details see Section~\ref{section::experiments}).

%%%%%%%%%%%%%%%%%%%%%%%%%%%%%%%%%%%%%%%%%%%%%%%%%%%%%%%%%%%%%%%%%%%%%%%%%%%%%%
\subsection{Contribution}

As already mentioned above, the major challenge of sparse feature extraction is to robustly identify a \emph{small} set of variables (non-zero components of $\omega$) that can be used to accurately classify unknown proteomics data (e.g. healthy or diseased) according to \eqref{eq:binarymodel}. % by learning from a given training set.
This paper introduces \emph{Sparse Proteomics Analysis} (\emph{SPA}), a novel framework for feature selection and classification. The key step of our method is based on \emph{1-bit compressed sensing} (cf. Section~\ref{section::CS}) and solves the following optimization problem:\footnote{Here, $\langle \cdot, \cdot \rangle$ again denotes the Euclidean scalar product.}
\begin{equation}\label{eq:PV-sparse}
\max_{\omega \in \R^d} \sum_{i=1}^n y_i\langle x_i, \omega \rangle\quad \text{subject to $\|\omega\|_1\le\sqrt{\lambda} $ and $\|\omega\|_2\le 1$,}
\end{equation}
where the regularization is now defined by two inequality constraints on the feature vector $\omega$.\footnote{For the sake of convenience, we formulate our algorithm as in \eqref{eq:PV-sparse}, but with some slight modifications, it could be equivalently stated in the form of \eqref{eqn:featuresel_intro}.} The above approach is motivated by the general theory of \emph{compressed sensing}, which was originally introduced by Donoho as well as by Cand\`es, Romberg, and Tao (cf. \cite{D, CT, CRT2}) and provides a modern framework for efficiently acquiring and processing high-dimensional (nearly) sparse signals (for more details see Section~\ref{section::CS}).

We shall verify the competitiveness of our method by applying it to several synthetic and real-world data-sets and comparing the results to those of other widely-used algorithms in this field. Although the core of the algorithm \eqref{eq:PV-sparse} is surprisingly simple, we will observe that SPA (including pre- and postprocessing steps) finds optimal feature vectors which are extremely sparse, allow for highly accurate classification, and are robust against noise.
In particular, for ``very-sparse'' situations, it even turns out that SPA outperforms the standard methods listed in Table \ref{tab::loss-and-reg-overview}.

Note that computational solutions to \eqref{eqn:featuresel_intro} or \eqref{eq:PV-sparse} are usually based on solving a convex program by standard optimization techniques, such as interior point methods. However, these methods sometimes scale poorly with increasing number of samples $n$ and data dimension $d$, as it is typically the case for -omics data analysis. Several strategies have been proposed in the literature to speed up the calculations, e.g., using stochastic decent (\cite{Genkin2007, Friedman2008, Efron2004, Koh2007, Wu2008}). In this article, we shall not focus on such computational issues but rather on providing a novel way of formalizing and solving the feature selection problem, namely in the context of compressed sensing.

Apart from the specific approach of \eqref{eq:PV-sparse}, it is a general concern of this work to promote the benefit of \emph{sparse} embedded methods.
In contrast to classical (univariate) approaches, such as statistical tests, the process of variable selection takes place in an automatic fashion here. In this way, a costly preprocessing (e.g., peak detection) as well as subsequent feature assessments can be avoided as much as possible.
Especially in a situation where only a very few samples are available, those additional steps may cause further instability and their success strongly relies on the specific data structure.
In fact, it was already succinctly emphasized by Vapnik in \cite[p.~12]{vapnik1998learning} that
\begin{center}
\emph{``If you possess a restricted amount of information for solving some problem, try to solve the problem directly and never solve the more general problem as an intermediate step. It is possible that the available information is sufficient for a direct solution but is insufficient for solving a more general intermediate problem.''}
\end{center}
This fundamental principle is precisely reflected by our viewpoint, which only makes a few (generic) assumptions on the underlying data model.
Finally, we would like to mention that recently, rigorous theoretical guarantees for sparse feature selection from MS data were shown in \cite{genzel2016selection}.
Using the novel idea of \emph{optimal problem representations}, the mathematical framework of \cite{genzel2016selection} even goes beyond the binary output scheme of \eqref{eq:binarymodel} and allows for a unified treatment of general observation and data models.

%%%%%%%%%%%%%%%%%%%%%%%%%%%%%%%%%%%%%%%%%%%%%%%%%%%%%%%%%%%%%%%%%%%%%%%%%%%%%%
\subsection{Outline of this Paper}
We start by shortly reviewing the background of \emph{compressed sensing} in Section \ref{section::CS}, and then describe our novel feature selection approach SPA in detail (Section \ref{section::spa-algo}). Finally, we present several benchmark results in Sections \ref{section::experiments} and \ref{section::experiments_real} for simulated and real data-sets and compare them to current state-of-the-art algorithms.

\section{Background: Compressed Sensing}
\label{section::CS}
%%%%%%%%%%%%%%%%%%%%%%%%%%%%%%%%%%%%%%%%%%%%%%%%%%%%%%%%%%%%%%%%%%%%%%%%%%%%%%
\subsection{Compressed Sensing-Based Data Analysis}
\label{subsec:cs}

In its most simple form, \emph{compressed sensing} (\emph{CS}) studies the recovery of an unknown vector $x \in \R^d$ from \emph{linear measurements} $y = Ax$.
Here, $A\in\R^{n\times d}$ is an $(n\times d)$-matrix and the entries of $y\in\R^n$ contain the measurements.
The major challenge is now to design the measurement process $A$ in such a way that the number of measurements $n$ is as small as possible and, at the same time, $x$ is still (uniquely) recoverable from $y$. Thus, we are asking for the maximal \emph{compressibility} of $x$ by linear measurements.

Obviously, when $n \ll d$, we require some additional information to obtain a unique solution of $y = Ax$. The prior information on $x$ which is studied in compressed sensing is the assumption of \emph{sparsity}, i.e., most coefficients of $x$ are assumed to be zero, or at least very small.
One naive approach to incorporate this additional property is to search for the sparsest solution of $Az = y$:\footnote{Here, $\|z\|_0 := \#\{ i \mid z_i \neq 0  \}$ simply counts the number of non-zero elements of $z = (z_1, \dots, z_d) \in \R^d$.}
\begin{equation}\label{eq:l0-intro}
\mathop{\rm arg min}_{z\in\R^d} \|z\|_0  \quad\text{subject to}\quad Az=y.
\end{equation}
Unfortunately, this problem is non-convex and cannot be efficiently solved in general. Therefore, one usually replaces \eqref{eq:l0-intro} by its \emph{convex relaxation}, which is also known as the \emph{basis pursuit} (\cite{CDS}):
\begin{equation}\label{eq:l1-intro}
\mathop{\rm arg min}_{z\in\R^d}\| z \|_1 \quad\text{subject to}\quad Az=y.
\end{equation}

One of the first key results in compressed sensing states that, if $A \in \R^{n\times d}$ is chosen \emph{randomly}, e.g., with independent and identically distributed Gaussian entries, and $n = O(s \cdot \log(d/s))$, then (with ``high probability'') every $s$-sparse vector $x$ (i.e., $\|x\|_0 \leq s$) can be uniquely recovered from \eqref{eq:l1-intro}.
The most surprising fact is that the number of required measurements $n = O(s \cdot \log(d/s))$ only logarithmically depends on the (possibly large) dimension $d$ of the ambient space.
% is almost of the order of the sparsity level $s$.
Hence, random measurement processes indeed allow for a very strong compression of sparse vectors (see also \cite{D, CT, CRT2} for more details).

In order to consider more complicated situations, the stability and robustness of the basis pursuit algorithm was extensively studied. Various theoretical results and
numerical experiments show that this algorithmic approach can also be applied for the stable recovery of
vectors which are only nearly sparse, as well as to noisy measurements of the form $y = Ax + \eta$.
To obtain a robust version of \eqref{eq:l1-intro}, one may replace its equality constraint by $\| Az - y \|_2 \leq \epsilon$ for some appropriate noise level $\epsilon > 0$. Not very surprisingly, this approach is also closely related to the Lasso introduced by Tibshirani in \cite{T} (see also \eqref{eqn:featuresel_intro} and Table~\ref{tab::loss-and-reg-overview}).

%%%%%%%%%%%%%%%%%%%%%%%%%%%%%%%%%%%%%%%%%%%%%%%%%%%%%%%%%%%%%%%%%%%%%%%%%%%%%%
\subsection{1-Bit Compressed Sensing}
%%%%%%%%%%%%%%%%%%%%%%%%%%%%%%%%%%%%%%%%%%%%%%%%%%%%%%%%%%%%%%%%%%%%%%%%%%%%%%
\label{subsec:1bitcs}

In many practical scenarios, especially when working with computers, there is no way to represent real numbers exactly. Thus, it is reasonable to assume that the measurement vector $Ax$ is acquired in a \emph{quantized} (and therefore non-linear) fashion. The most extreme form directly leads to \emph{1-bit measurements}, i.e., only the signs of $Ax$ are known:
\begin{equation}\label{eq:1bitmeas}
	y_i={\sign}(\langle a_i, x \rangle),\quad \quad i=1,\dots,n,
\end{equation}
where $a_1,\dots,a_n\in\R^d$ are the rows of the measurement matrix $A \in \R^{n \times d}$. As in classical compressed sensing, we are asking for an appropriate recovery of $x$ from \eqref{eq:1bitmeas} using as few measurements as possible. This challenge was originally considered in \cite{BB} as \emph{1-bit compressed sensing}, and has been extensively studied in \cite{PV1, PV2}.

A surprisingly simple convex recovery approach was proposed by Plan and Vershynin in \cite{PV2}:
\begin{equation}\label{eq:PV}
\max_{z \in \R^d} \sum_{i=1}^n y_i \langle a_i, z \rangle\quad \text{subject to $\|z\|_1\le\sqrt{\lambda}$ and $\|z\|_2\le 1$,}
\end{equation}
where $\lambda > 0$ denotes the sparsity-controlling parameter. To get some intuition, we first note that we have $y_i={\sign}(\langle a_i, x \rangle)$ if and only if $y_i \langle a_i, x \rangle > 0$ holds. Hence, maximizing the sum in \eqref{eq:PV} will ensure the consistency of many measurements $i \in \{1, \dots, n\}$, according to \eqref{eq:1bitmeas}. However, the total consistency is not enforced so that \eqref{eq:PV} indeed allows for noisy inputs $y$ that do not satisfy \eqref{eq:1bitmeas}. On the other hand, the constraint of \eqref{eq:PV} promotes sparsity of the final outcome. To see this, we may consider the set $S_{d,\lambda} := \{ z \in \R^d : \|z\|_0\le \lambda, \|z\|_2 \leq 1 \}$ and observe that (cf. \cite[Sec.~III]{PV2})\footnote{Here, $\operatorname{conv}(S)$ denotes the convex hull of the set $S \subset \R^d$.}
\begin{equation*}
\operatorname{conv}(S_{d,\lambda}) \subset \{ z \in \R^d : \|z\|_1\le\sqrt{\lambda}, \|z\|_2\le 1 \} \subset 2 \operatorname{conv}(S_{d,\lambda}).
\end{equation*}
This means that \eqref{eq:PV} optimizes over a convex relaxation of the set $S_{d,\lambda}$ which contains all $\lambda$-sparse vectors in the unit ball. For more details, see also \cite{PV2}.
The main statement of \cite{PV2} proves that the robust 1-bit compressed sensing algorithm \eqref{eq:PV} indeed allows for an appropriate recovery of sparse vectors, using only $n = O(\lambda \cdot \log(d/\lambda))$ measurements. Moreover, it is surprisingly robust against several types of noise, including (random) bit-flips of the labels.

\begin{remark*}
The minimized functional of \eqref{eq:PV} is closely related to the hinge loss which is used for SVMs (cf. Table~\ref{tab::loss-and-reg-overview}). Indeed, without rejecting the negative part of the hinge loss, we would precisely end up with the objective functional in \eqref{eq:PV}.

The constraint of \eqref{eq:PV}, on the other hand, can be regarded as a combined $\ell_1$-$\ell_2$-condition, where the tuning parameter $\lambda$ controls the desired level of sparsity of the minimizer. This type of regularization strongly resembles the idea of \emph{elastic nets}, originally proposed by Zou and Hastie in \cite{elastic}.
\end{remark*}

%%%%%%%%%%%%%%%%%%%%%%%%%%%%%%%%%%%%%%%%%%%%%%%%%%%%%%%%%%%%%%%%%%%%%%%%%%%%%%
\subsection{Why Compressed Sensing?}
%%%%%%%%%%%%%%%%%%%%%%%%%%%%%%%%%%%%%%%%%%%%%%%%%%%%%%%%%%%%%%%%%%%%%%%%%%%%%%
\label{subsec:csml}

At a first sight, the main challenges of compressed sensing and machine learning (ML) seem to be very different. In compressed sensing, we intend to design a measurement process $A$ in order to \emph{compress} a vector $x$, whereas in machine learning, the training data is already contained in the rows of $A$ and we are rather willing to \emph{explain} the observations $y$ by some appropriate vector $x$. However, in both areas we are asking for a (sparse) recovery from a certain type of measurement. Indeed, a \emph{linear regression} in ML exactly corresponds to classical CS model (see Subsection~\ref{subsec:cs}), and a \emph{classification} problem is actually equivalent to 1-bit CS (see Subsection~\ref{subsec:1bitcs}).

Therefore, it is not very surprising that the applied algorithms for compressed sensing and machine learning resemble each other, and that theoretical results in both fields rely on the same mathematical foundations (concentration of measure, convex geometry, etc.). Unfortunately, both communities only rarely interacted with each other.
In this paper, we would like to emphasize the viewpoint of compressed sensing, in particular, because it is still not very common for the classification tasks that we deal with.

With the recent progress in compressed sensing and related areas as low-rank matrix recovery or quantized CS, also new algorithms like nuclear norm minimization or 1-bit CS have been proposed.
Although these methods are typically motivated by theoretical studies, they perform also very well for real-world data.
In general, we believe that these alternative perspectives allow for deeper theoretical insights, finally leading to the improvement of the classical ($\ell_1$-based) tools from machine learning.

For an extensive introduction to compressed sensing, we refer to \cite{DDEK, FR}. As we already mentioned above, comparing this text to literature from statistical learning theory (see \cite{Buhlmann2011} for example), the reader will quickly notice many interesting connections between both fields.

%%%%%%%%%%%%%%%%%%%%%%%%%%%%%%%%%%%%%%%%%%%%%%%%%%%%%%%%%%%%%%%%%%%%%%%%%%%%%%
%\clearpage
% \section{Sparse Proteomics Analysis (\methodname): A new method for feature detection in high-dim Proteomics data}
\section{Sparse Proteomics Analysis (SPA)}
\label{section::spa-algo}
%%%%%%%%%%%%%%%%%%%%%%%%%%%%%%%%%%%%%%%%%%%%%%%%%%%%%%%%%%%%%%%%%%%%%%%%%%%%%%
%%%%%%%%%%%%%%%%%%%%%%%%%%%%%%%%%%%%%%%%%%%%%%%%%%%%%%%%%%%%%%%%%%%%%%%%%%%%%%

In this section, we present the details of our novel framework which is based on the ideas of 1-bit compressed sensing presented in the previous section. The first part provides a mathematical formulation of the feature selection problem as well as a brief overview of the steps that are performed in SPA. The rest of this section is then devoted to a detailed description and discussion of the single steps.

%\begin{wrapfigure}{r}{0.45\textwidth}
%\vspace{-0.9cm}
%\fbox{\parbox{0.43\textwidth}{
%\begin{center}
%\includegraphics[width=5cm]{figures/toy_spectra.png}
%\caption{A simple toy example with three peaks. The plot shows data of two classes indicated as red and blue curves. Our goal is to %identify a small set of positions that allows for a robust distinction of the two classes. Here, the significant points are clearly $\{5, 10, %15\}$.}
%\label{fig:simple_example_three_peaks}
%\end{center}
%}}
%\end{wrapfigure}

%%%%%%%%%%%%%%%%%%%%%%%%%%%%%%%%%%%%%%%%%%%%%%%%%%%%%%%%%%%%%%%%%%%%%%%%%%%%%%
\subsection{Setting and Overview}
%%%%%%%%%%%%%%%%%%%%%%%%%%%%%%%%%%%%%%%%%%%%%%%%%%%%%%%%%%%%%%%%%%%%%%%%%%%%%%

As already mentioned in the introduction, we assume that our learning process is \emph{supervised}, i.e.,  we know which spectrum belongs to the class of healthy ($y_i=+1$) and diseased ($y_i=-1$) samples in advance. If the data vectors  $x_i \in \R^d$, $i=1,\ldots,n$ are 
mass spectra, the indices $j = 1, \dots, d$ of $x_i = (x_{i,1}, \dots, x_{i,d})$ correspond to the $m/z$-values\footnote{$m/z$ is the unit for the mass-to-charge ratio.} and its entries represent the intensities.
The non-zero entries of the feature vector $\omega_0 = (\omega_{0,1}, \dots, \omega_{0,d})\in \R^d$ describe the location of the disease fingerprints and its respective values the significance of these features.

In the setting of classical learning theory, we are asking for a hyperplane $\{\omega_0\}^\perp$ which correctly separates most of the data points $x_i$ labeled by $y_i$. More precisely, this means\footnote{Compared to Section \ref{section::CS}, we are now using the standard notations from learning theory. In particular, the measurement vectors are denoted by $x_i$ (instead of $a_i$) and the feature vector is $\omega_0$ (instead of $x$).}
\begin{equation}
\label{eq:measurement}
y_i = \sign(\langle x_i, \omega_0 \rangle)\quad \text{for ``many'' samples $i = 1, \dots, n$.}
\end{equation}
Equivalently, we can view \eqref{eq:measurement} as a problem from 1-bit compressed sensing (cf. Section~\ref{subsec:csml}), i.e., we have acquired noisy 1-bit measurements and are now looking for a sparse recovery.

In the development of SPA, we have primarily focused on the latter perspective,
and therefore, the 1-bit recovery program \eqref{eq:PV} forms the key step of our algorithm:

% \fbox{\parbox{.9\textwidth}{
% \begin{algo}[SPA Overview]\leavevmode
\algobox[\textbf{SPA at a Glance}]{.95\textwidth}{\label{algo:spa}

\emph{Input:} Raw data samples $\{(x_i,y_i)\}_{i=1,\dots,n}$

\emph{Output:} Sparse feature vector $\tilde{\omega} \in \R^d$

\textbf{Preprocessing:}
\begin{enumerate}
\item[1:]
	Normalize data to make the spectra comparable.
\item[2:]
	Perform smoothing by a convolution with Gaussian density.
\item[3:]
	Standardize data.
\end{enumerate}
% The $x_i$ are the spectra resulting from these preprocessing steps.\\
\textbf{Sparse Feature Selection:}
\begin{enumerate}
\item[4:]
	Perform 1-bit CS optimization \eqref{eq:PV} to find feature vector $\hat\omega$.
\end{enumerate}
\textbf{Postprocessing:}
\begin{enumerate}
\item[5:]
	Detect the connected components of $\hat\omega$ to obtain a sparsified version $\tilde\omega$.
\item[6:]
	(Optional) Reduce dimension by projecting data onto the feature space.
\end{enumerate}}
% \end{algo}\vspace{-1\baselineskip}}}
% \vspace{\baselineskip}

%%%%%%%%%%%%%%%%%%%%%%%%%%%%%%%%%%%%%%%%%%%%%%%%%%%%%%%%%%%%%%%%%%%%%%%%%%%%%%
\subsection{Algorithmic Details}
\label{section::AlgorithmicDetails}
%%%%%%%%%%%%%%%%%%%%%%%%%%%%%%%%%%%%%%%%%%%%%%%%%%%%%%%%%%%%%%%%%%%%%%%%%%%%%%

In the following, we are going to specify and discuss the single steps of Algorithm~\ref{algo:spa}.

%%%%%%%%%%%%%%%%%%%%%%%%%%%%%%%%%%%%%%%%%%%%%%%%%%%%%%%%%%%%%%%%%%%%%%%%%%%%%%
\subsubsection*{Step 1: Normalization of the Data}

This step heavily depends on the underlying acquisition method of the data. Every spectrum $x_i \in \R^d$ is normalized by a certain scaling factor $\mu_i > 0$, i.e., $x_i \mapsto \mu_i x_i$ for $i = 1, \dots, n$. The individual scalars $\mu_i$ should be chosen such that the resulting data vectors are ``comparable.''

For example, when we assume that the data are acquired by MALDI-TOF-MS as described in Fig.~\ref{fig1}, it seems to be quite natural to normalize them by the total ion count. Mathematically, this means that we would divide every spectrum by its $\ell_1$-norm, i.e., we choose $\mu_i = 1 / \|x_i\|_{1}$.\

%%%%%%%%%%%%%%%%%%%%%%%%%%%%%%%%%%%%%%%%%%%%%%%%%%%%%%%%%%%%%%%%%%%%%%%%%%%%%%
\subsubsection*{Step 2: Smoothing by Gaussian Density}

We already pointed out that one major challenge is the strong noise within the raw data. Therefore, it is crucial to perform some noise reduction before trying to extract features.
For this purpose, we suggest a simple smoothing strategy by a Gaussian density:

Let $G_\sigma$ denote the (centered) \emph{Gaussian density function} with fixed standard deviation $\sigma > 0$, that is,
\begin{equation*}
G_\sigma (t) = \frac{1}{\sqrt{2\pi \sigma^2}} \exp\left( -\frac{t^2}{2\sigma^2} \right), \quad t \in \R.
\end{equation*}
The smoothed spectra $\tilde{x}_i = (\tilde{x}_{i,1}, \dots, \tilde{x}_{i,d}) \in \R^d$ are then obtained by a discrete convolution
\begin{equation}\label{eq:gaussdic}
\tilde{x}_{i,k} := (x_i \ast G_\sigma)_k = \sum_{l = 1}^d x_{i,l} \cdot G_\sigma(k - l), \quad k = 1, \dots, d, \quad i = 1, \dots, n.
\end{equation}
Using the fast Fourier transform (FFT), this computation can be performed quickly with $O(n d \log(d))$ operations.
In a very simplified scenario, a spectrum can be written as the sum of Gaussian-shaped peaks plus some baseline noise in each mass channel. Since the convolution of two Gaussian densities is again Gaussian, the original (local) structure of the spectra is essentially preserved in $\tilde{x}_i$, whereas the noise of $x_i$ is significantly reduced. Note that the deviation $\sigma > 0$ serves as a tuning parameter of the algorithm. A good choice of $\sigma$ clearly depends on the nature of the data; usually it depends on the noise level as well as on the (average) width of the peaks.

Finally, we would like to emphasize another interesting interpretation of the above smoothing approach: The convolution in \eqref{eq:gaussdic} can be written as a scalar product of $x_i$ with the shifted Gaussian density $G_\sigma(\cdot - k)$ (note that $G_\sigma$ is symmetric), that is, $\tilde{x}_{i,k} = \langle x_i, G_\sigma(\cdot - k) \rangle$. Thus, the entries of $\tilde{x}_i$ are actually the \emph{analysis coefficients} of the \emph{Gaussian dictionary} $\{G_\sigma(\cdot - k) \mid k = 1, \dots, d \}$. The perspective of analyzing data by a \emph{dictionary} offers several opportunities for generalization. For instance, one could also consider (redundant) dictionaries with more than one standard deviation or more sophisticated functions than $G_\sigma$.

%%%%%%%%%%%%%%%%%%%%%%%%%%%%%%%%%%%%%%%%%%%%%%%%%%%%%%%%%%%%%%%%%%%%%%%%%%%%%%
\subsubsection*{Step 3: Standardizing the Data}

The 1-bit optimization of \eqref{eq:PV} does not incorporate a bias term. Hence, it is necessary to center the data first. For this, we compute the \emph{mean spectrum}\footnote{Actually, we use the smoothed data vectors $\tilde{x}_i$ from Step 2 as input for this computation. But in order to keep the notation simple, we still write $x_i$. This convention holds also for all forthcoming steps.}
\begin{equation*}
\bar x := \tfrac{1}{n} \sum_{i = 1}^n x_i \in \R^d,
\end{equation*}
i.e., $\bar{x}_k$ contains the average of the $k$-th entry of all spectra. %Then, we obtain the \emph{centered  spectra} by subtracting the average:
%\begin{equation*}
%\bar{x}_i := x_i - \bar x, \quad i = 1, \dots, n.
%\end{equation*}
The spectra are further scaled by dividing the non-constant features by their \emph{standard deviation}
\begin{equation*}
	\sigma_j := \sqrt{\tfrac{1}{n}\sum_{i = 1}^{n}\left( x_{i,j} - \bar{x}_j\right)^2},\quad j = 1, \dots, d.
\end{equation*}
The \emph{standardized spectra} $\check{x}_i = (\check{x}_{i,1}, \dots, \check{x}_{i,d}) \in \R^d$ are then obtained by
\begin{equation*}
	\check{x}_{i,j}:=\frac{x_{i,j} - \bar{x}_j}{\sigma_j},\quad i = 1, \dots , n, \quad  j = 1, \dots , d.
\end{equation*}
In this way, all feature variables are centered and have an empirical standard deviation equal to $1$, so that they get equally weighted in the selection process.

%%%%%%%%%%%%%%%%%%%%%%%%%%%%%%%%%%%%%%%%%%%%%%%%%%%%%%%%%%%%%%%%%%%%%%%%%%%%%%
\subsubsection*{Step 4: Sparse Feature Selection}

We are now ready to perform the actual feature extraction step, using the 1-bit recovery method presented in Subsection~\ref{subsec:1bitcs}:

\algobox[1-Bit Compressed Sensing]{.95\textwidth}{\label{algo:onebitfeature}

	\emph{Input:} Samples $\{(x_i,y_i)\}_{i=1,\dots,n}$, sparsity parameter $\lambda > 0$, threshold $\epsilon > 0$

	\emph{Output:} Estimated feature vector $\hat\omega = (\hat{\omega}_1, \dots, \hat{\omega}_d) \in \R^d$

	\emph{Compute:}
	\begin{align}\label{eq:onebitfeature}
		\text{1:} \quad & \hat{\omega}' = \argmax_{\omega \in \R^d} \sum_{i = 1}^n y_i \langle x_i, \omega \rangle \quad \text{subject to $\|\omega\|_1 \le \sqrt{\lambda}$ and $\|\omega\|_2 \le 1$.} \\
		\label{eq:hardthres}
		\text{2:} \quad & \hat{\omega}_k = \begin{cases} \hat{\omega}'_k, & \text{if } \abs{\hat{\omega}'_k} > \epsilon, \\
			0, & \text{otherwise,}
		\end{cases} \quad k = 1, \dots, d.
	\end{align}
}

The second part (in \eqref{eq:hardthres}) is a simple hard thresholding that tries to eliminate computational inaccuracies by setting almost zero entries of $\hat{\omega}'$ to $0$ ($\epsilon$ is usually very small, e.g.,  $\sim 10^{-3}$).

The actual feature selection takes place in \eqref{eq:onebitfeature}. Recalling the observation model from \eqref{eq:measurement}, we conclude that the $i$-th sample is correctly classified by a vector $\omega$ if and only if $y_i \langle x_i, \omega \rangle > 0$.
Hence, the objective functional of \eqref{eq:onebitfeature} will be particularly large if sufficiently many samples are correctly classified by $\omega$.
However, a consistent prediction of \emph{all} measurements (i.e., $y_i = \sign(\langle x_i, \omega\rangle)$ for all $i = 1, \dots, n$) is not strictly enforced, and therefore, our strategy enjoys a certain robustness against (random) perturbation of the model \eqref{eq:measurement}. This could occur in practice, for example, when a training sample was wrongly classified from the very beginning. On the other hand, the constraint of \eqref{eq:onebitfeature} guarantees that the maximizer will be ``effectively'' sparse (depending on the choice of the sparsity parameter $\lambda > 0$). This intuition indicates that the estimator $\hat{\omega}$ will be indeed a sparse vector allowing for an appropriate separation of the two classes.

%%%%%%%%%%%%%%%%%%%%%%%%%%%%%%%%%%%%%%%%%%%%%%%%%%%%%%%%%%%%%%%%%%%%%%%%%%%%%%
\subsubsection*{Step 5: Detecting the Connected Components}

One advantage of Algorithm~\ref{algo:onebitfeature} is that it does not make any assumptions on the structure of the data vectors $x_i$. Hence, it might be even suited for much more general types of data. However, its ``universality'' comes with the drawback that the characteristic peak structure of MS data is not captured at all. In fact, a spectrum does not consist of sharp spikes but rather wide-spread Gaussian shaped peaks. Hence, if the algorithm finds a significant feature position, say at the maximum of some peak, it usually tends to select also those features which are close to this position. Such a behavior is not very surprising, because nearby features are highly correlated to the maximum of the peak, and therefore, they may allow for a good separation as well.

\begin{figure}
% peak_connected_interval.png
\centering
\includegraphics[width=4cm]{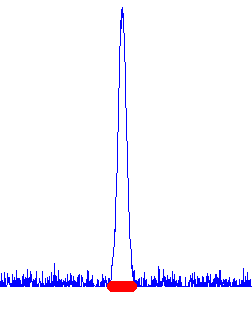}
\caption{The red stripe indicates the support of $\hat{\omega}$. Relevant features usually occur as intervals and not as isolated points.}
\label{fig3}
% fig:peak_connected_interval
\end{figure}

Empirical results have shown that this process of selection ``evolves'' in a continuous fashion when changing the sparsity level $\lambda$. As a consequence, the support of a feature vector $\hat{\omega}$ from Algorithm~\ref{algo:onebitfeature} typically consists of a few connected ``intervals'' (consecutive sequences of indices) which are centered around the selected peaks (see also Fig.~\ref{fig3}). The actual sparsity of $\hat{\omega}$ should be therefore measured by means of its connected intervals and not by simply counting its non-zero entries.

For this reason, we may easily improve the sparsity of $\hat{\omega}$ by reducing every interval to its most significant entry:\footnote{Here $\supp(\hat{\omega}) = \{ k \mid \hat{\omega}_k \neq 0\}$ denotes the support of $\hat{\omega}$, i.e., the set of indices corresponding to its non-zero entries.}

\algobox[Sparsification of $\hat{\omega}$]{.95\textwidth}{

	\emph{Input:} (Sparse) feature vector $\hat{\omega} = (\hat{\omega}_1, \dots, \hat{\omega}_d) \in \R^d$

	\emph{Output:} Sparsified version $\tilde\omega = (\tilde{\omega}_1, \dots, \tilde{\omega}_d) \in \R^d$

	\emph{Compute:}
	\begin{enumerate}
	\item[1:]
		Find the connected components $A_1, \dots, A_N \subset \supp(\hat\omega)$ of $\hat\omega$.
	\item[2:]
		For every $l = 1, \dots, N$ do the following:
		
		Set all entries of $\hat{\omega}$ in $A_l$ to $0$, except from $\argmax_{k \in A_l} |\hat{\omega}_k|$.
	\item[3:]
		The resulting vector is $\tilde\omega$.
	\end{enumerate}
}

%%%%%%%%%%%%%%%%%%%%%%%%%%%%%%%%%%%%%%%%%%%%%%%%%%%%%%%%%%%%%%%%%%%%%%%%%%%%%%
\subsubsection*{Step 6: Dimension Reduction}
\label{subsec:dim-reduction}

This final (optional) step does not involve any further computations but shows how to proceed with our result $\tilde\omega$.
As mentioned before, the main purpose of SPA is not just to classify (unknown) samples, but rather to reduce the data to its significant entries (dimensions).
Indeed, we may use $\tilde\omega$ for a \emph{dimension reduction}: Let $x = (x_1, \dots, x_d) \in \R^d$ be some (possibly unknown) data vector.
Then, we can project $x$ onto the selected feature positions of $\supp(\tilde\omega)$. More precisely, all entries that do not belong to $\supp(\tilde\omega)$ are set to $0$:
\begin{equation}
\hat{x}_k := \begin{cases} x_k, & k \in \supp(\tilde\omega), \\
	0, & \text{otherwise,}
\end{cases}\quad k = 1, \dots, d.
\end{equation}
The resulting data vector $\hat{x} = (\hat{x}_1, \dots, \hat{x}_d) \in \R^d$ is now trivially embedded into a low-dimensional space of dimension $\# \supp(\tilde\omega)$.\footnote{In practice, one would simply reject all indices that are not contained in $\supp(\tilde\omega)$.} But it still contains the most important information which has been found by the above algorithm.
Note that we have not made any use of the actual values of $\tilde\omega$ but merely of its support.

By this projection, we may reduce the danger of overfitting. In particular, by working in a low-dimensional space, a large tool set from \emph{machine learning} is now available for classification and clustering. But how to explicitly proceed with the data heavily depends on the specific application and is therefore not part of SPA.

%%%%%%%%%%%%%%%%%%%%%%%%%%%%%%%%%%%%%%%%%%%%%%%%%%%%%%%%%%%%%%%%%%%%%%%%%%%%%%
% \clearpage
\section{Experimental Results: Feature Selection from Simulated Data-Sets}
\label{section::experiments}

In this section, we assess our framework of SPA with regard to a typical situation in mass-spectrometry analysis: We would like to extract discriminating features from MS data with respect to two groups (e.g., healthy and diseased patients). A major difficulty is usually that only a small number of measurements (observations) is available. Building on this, we ask for the following: Given a simulated data-set for which the position and number of discriminating peaks are known (this will be called $\omega_0$ below), how many samples are needed to identify these features with high accuracy?

We shall compare our results to the widely used state-of-the-art algorithms LIBLINEAR ($\ell_1$-regularized SVM) and the standard MATLAB implementation of Lasso.

\subsection{Creating a Simulated Data-Set}

\begin{figure}
	
	\captionsetup[subfigure]{justification=centering}	
	\centering
	\begin{subfigure}[b]{0.32\textwidth}
		\includegraphics[scale = 0.32]{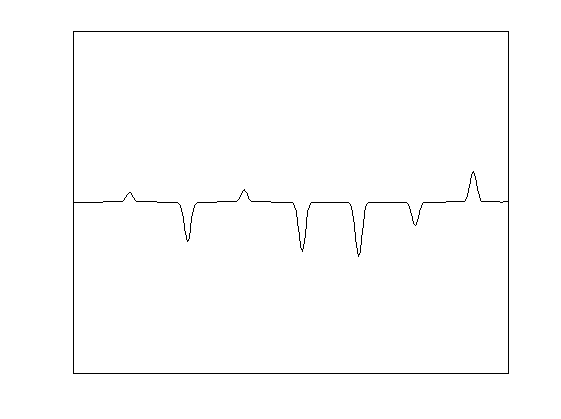}
		\subcaption{}

	\end{subfigure}
	\begin{subfigure}[b]{0.32\textwidth}
		\includegraphics[scale = 0.32]{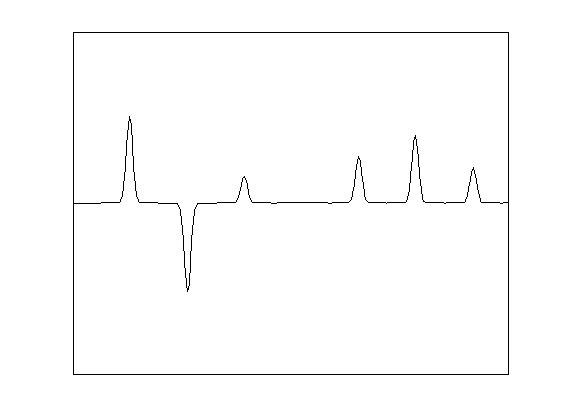}
		\subcaption{}
		
	\end{subfigure}
	\begin{subfigure}[b]{0.32\textwidth}
		\includegraphics[scale = 0.32]{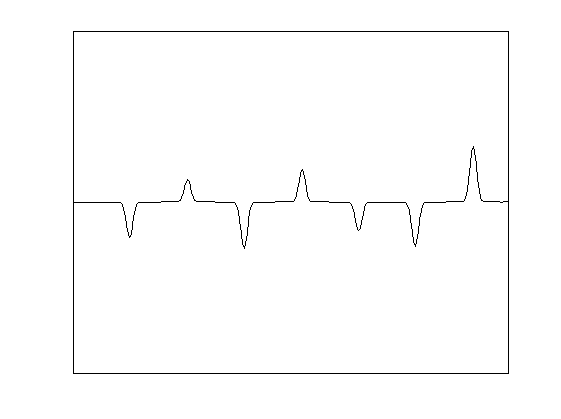}
		\subcaption{}
		
	\end{subfigure}	
	
	\begin{subfigure}[b]{0.32\textwidth}
			\includegraphics[scale = 0.32]{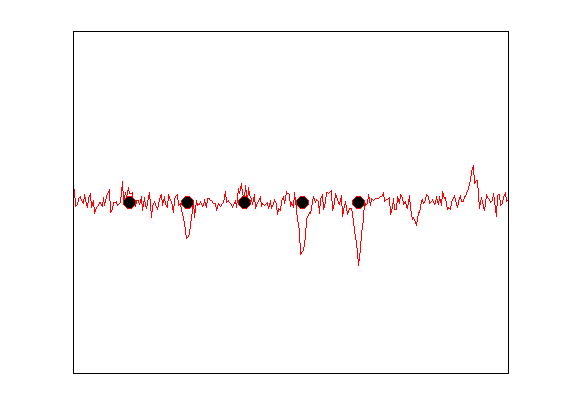}
			\subcaption{}
			
	\end{subfigure}	
	\begin{subfigure}[b]{0.32\textwidth}
			\includegraphics[scale = 0.32]{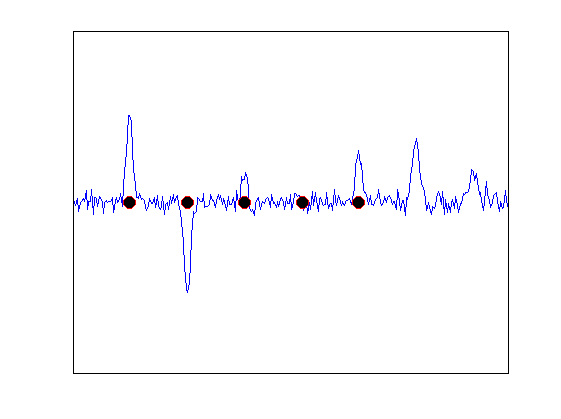}
			\subcaption{}
			
	\end{subfigure}
	\begin{subfigure}[b]{0.32\textwidth}
			\includegraphics[scale = 0.32]{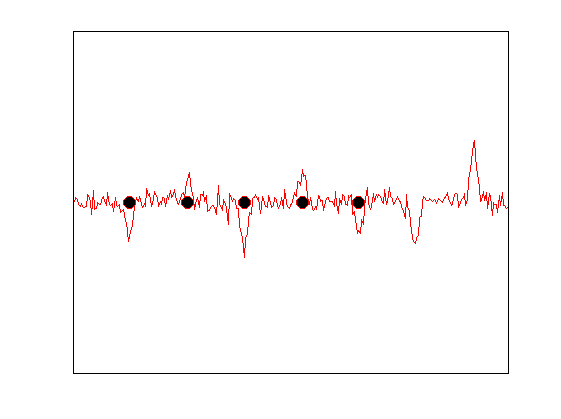}
			\subcaption{}
			
	\end{subfigure}
	\caption{Illustration of the generated data instances. (a)--(c): First seven equidistant Gaussian peaks that are located in fixed positions in each of the three data instances; (d)--(f): Visualization of the data instances from (a)--(c) with additive noise with standard deviation $\sigma = 0.1$, where the positions of the five condition positive peaks are highlighted by black dots. The blue and red colors indicate the different classes which are determined by the observation process of \eqref{eq:simdata1bitmeas}.}
	\label{data set}
\end{figure}

We assume that our sample set $\{(x_i, y_i)\}_{i = 1, \dots, n} \subset \R^d \times \{-1, +1\}$ follows a certain joint random distribution $(X,Y)$, where each sample is independently drawn.
In order to make the problem tractable, let us make two model assumptions on $X$ and $Y$. First, the mass spectra $X$ are generated as follows:
$$
	x_i = \sum_{m = 1}^{M} s_i^m a^m + n_i, \quad i = 1, \dots, n,
$$
where $s_i^m \in \R^d$ determines the (random) amplitude of the $m$-th peak, $a^m \in \R^d$ specifies its position and shape, and $n_i \in \R$ represents the low-amplitude baseline noise.
We shall assume that the amplitudes and the noise are Gaussian, that is, $s_i := (s_i^1,...,s_i^M) \sim \mathcal{N}(0,\Sigma)$ with $\Sigma \in \R^{M \times M}$ positive definite and $n_i \sim \mathcal{N}(0,\sigma^2I)$ with $\sigma > 0$.
Note that the generated data might have negative components. This does not mimic the structure of real-world mass spectra which is always non-negative.
However, since centering is part of our preprocessing anyway (cf. Step 3 in Subsection~\ref{section::AlgorithmicDetails}), the assumption of mean-zero amplitudes is quite natural. The (disease) labels $Y$ are then simply modeled as 1-bit observations (see also \eqref{eq:measurement})
\begin{equation}\label{eq:simdata1bitmeas}
	y_i = \sign ( \langle x_i, \omega_0 \rangle), \quad i = 1, \dots, n,
\end{equation}
where $\omega_0 \in \R^d$ is the sparse ground-truth feature vector, which we intend to estimate.
In the following, each non-zero entry of $\omega_0$ is located at the center of a specific peak (see Fig.~\ref{data set}(d)--(f)), so that $\supp(\omega_0)$ actually determines all biologically relevant peaks (molecular structures).
Since $\Sigma$ is invertible (i.e., the features are linearly independent), this collection of peaks is an optimal fingerprint in the sense that removing or adding any feature variable would decrease the prediction accuracy (with respect to the ``perfect'' model of \eqref{eq:simdata1bitmeas}).
% , which is supported within the few peaks. Since $\Sigma$ is invertible (i.e., no feature can be linearly combined by the others) the support of $\omega_0$ is unique. Thus, there exists no subset which gives the same classification performance. If we would remove a feature from $\omega_0$, the prediction accuracy would indeed decrease.

In our experiments, we create data-sets $ x_1, \dots, x_n \in \R^{8192} $, each one consisting of $200$ equidistant peaks (atoms $a^m$) shaped like Gaussian density function of width $10$.
The vector $\omega_0 \in \R^{8192}$ is chosen to have five non-zero components, which means that only five prechosen peaks were used to generate the labels $y_1, \dots, y_n$. Hereafter, we will refer to these as \emph{condition positive peaks}. Fig.~\ref{data set} shows three different data instances magnifying only the first seven peaks, generated in the described way.
In order to verify our method, we will use two types of data-sets DS1 and DS2 which only differ in their correlation matrix $\Sigma$.
For DS1, $\Sigma$ is chosen to be the identity matrix. This implies that the heights of all of the 200 peaks are standard Gaussian random variables.
For DS2, we have chosen three pairs of negative peaks to be positively correlated and in addition, one condition positive peak was chosen to be positively correlated with one of the negative peaks. Thus, there are a few entries of value $0.8$ off the main diagonal in $\Sigma$.
To test the algorithm's performance increasing amount of Gaussian noise $n_i \sim \mathcal{N}(0,\sigma^2)$  with $\sigma = \{0.1,0.3\}$ was added to DS1 and DS2. These corresponds to signal-to-noise (SNR) ratio of 10, 3.33 repectively \footnote{Signal-to-noise ratio was calculated as $SNR = \frac{\operatorname{power\ of\ signal}}{\operatorname{power\ of\ noise}}$.}. The values of SNR are chosen to represent the behaviour of the algorithm up to the levels of noise that are normally found in MS data.

\subsection{Setup and Evaluation Criteria}

Let us recall the essential question of our experiments: Can we recover the support of $\omega_0$, and if so, how many samples do we need for that?
% We will measure the algorithm's capability of selecting an optimal feature set as well as its dependence on the number of observations in the data set.
For this purpose, we shall successively increase the number of available samples in the (training) data-set and examine whether SPA (or Lasso, or $\ell_1$-SVM) succeeds in recovering $\supp(\omega_0)$.
Since each of the considered algorithms involves a variable parameter, we have decided to perform an adaptive tuning for each problem instance. In fact, the sparsity parameter was chosen such that the resulting classifier $\tilde{\omega}$ matches the sparsity level of $\omega_0$.
But of course, this does not automatically imply that the supports of $\tilde\omega$ and $\omega_0$ completely coincide.\footnote{Due to the redundancy of the peak-associated feature variables (cf. Step 5 in Subsection~\ref{section::AlgorithmicDetails}), an estimated feature vector is considered to be equal to the ground-truth vector with some tolerance, which particularly depends on the width of the peaks.}
For each problem instance, the smallest sparsity parameter which resulted in a classifier with five non-zero entries was chosen in the following way: The initial value of the sparsity parameter for SPA and $\ell_1$-SVM (Lasso) was set to the value which corresponds to the classifier with less than (more than) five non-zero values\footnote{This difference arises from the implementation of Lasso.}. For SPA and $\ell_1$-SVM (Lasso), the sparsity parameter was increased for a preset step size until the outcome had five or more (five or fewer) non-zero entries. If the previous step provided a sparse classifier with strictly more than (strictly less than) five non-zero entries, the bisection method was used on the interval between the two last sparsity parameter values. The bisection method was used until the optimal sparsity parameter was found or the difference between the two consecutive parameters became smaller than a preset tolerance.

We will use a measure based on \emph{sensitivity}. Sensitivity, defined as\footnote{TP - true positives, i.e. correctly identified peaks\\ FP - false positives, i.e. incorrectly identified peaks\\ TN - true negatives, i.e. correctly rejected peaks\\ FN - false negatives, i.e. incorrectly rejected peaks}
\begin{equation*}
\operatorname{sens} := \tfrac{TP}{TP + FN}
\end{equation*}

is an appropriate measure for our objectives because it represents an algorithm's ability to detect the relevant features. Note that ideally, the number of condition positives ($TP+FN$) is equal to predicted condition positives ($TP+FP$). In such a situation, the \emph{precision}, given by $p := TP/(TP + FP)$ is equal to the sensitivity. However, in the presence of noise it is possible that the final selection encompasses several features which are associated with a single peak. This could lead to a precision value equal to $1$ if all of the selected values are declared as true positives, though some other true features remain undetected. Since for us, it is equally important to penalize both false positives and false negatives, we have chosen the sensitivity to be the main point of reference. A measure of similar importance is the \emph{specificity}, which is defined by
\begin{equation*}
\operatorname{spec} := \tfrac{TN}{FP + TN}.
\end{equation*}
Finally, due to the possibly imbalanced number of relevant features, we shall also take into account the so-called \emph{balanced accuracy}
\begin{equation*}
\operatorname{bacc} := \tfrac{\operatorname{sens} + \operatorname{spec}}{2}.
\end{equation*}

\subsection{Results}
Data-sets of sample sizes between $50$ and $350$ were generated as described above and each of the methods was performed for standardized input data. Note that the hard thresholding step described in \eqref{eq:hardthres} was also applied to the classifiers obtained from Lasso or $\ell_1$-SVM. Otherwise, any computational inaccuracy would completely destroy the sparsity structure of the results.

\begin{figure}
	\centering
		\includegraphics[scale = 0.4]{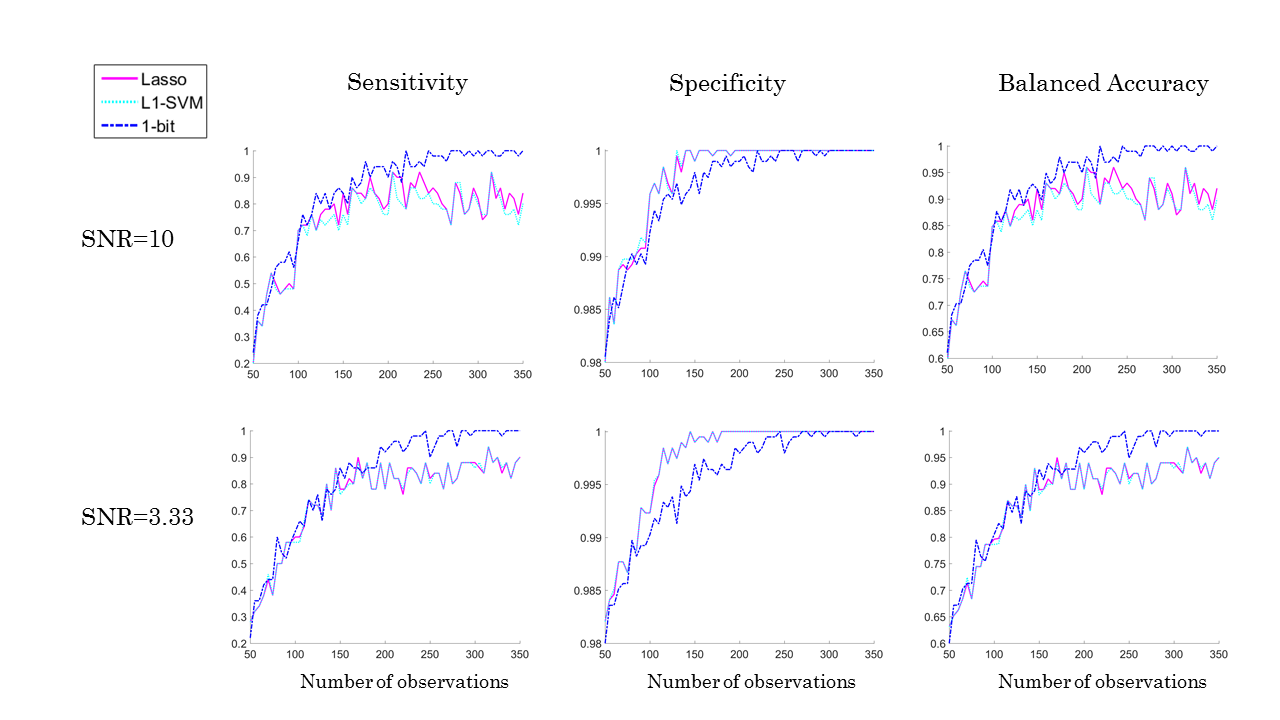}
	\caption{Comparison of numerical results for SPA ($= \text{$1$-bit CS}$), Lasso, and $\ell_1$-SVM on the data-set DS1 with SNR = 10, and 3.33, showed in the respective row. Note that the data consist of $5$ condition positive and $195$ condition negative peaks which are equidistantly located in the spectra.}
	\label{fig::DS1}
\end{figure}

\begin{figure}
	
		\includegraphics[scale = 0.40]{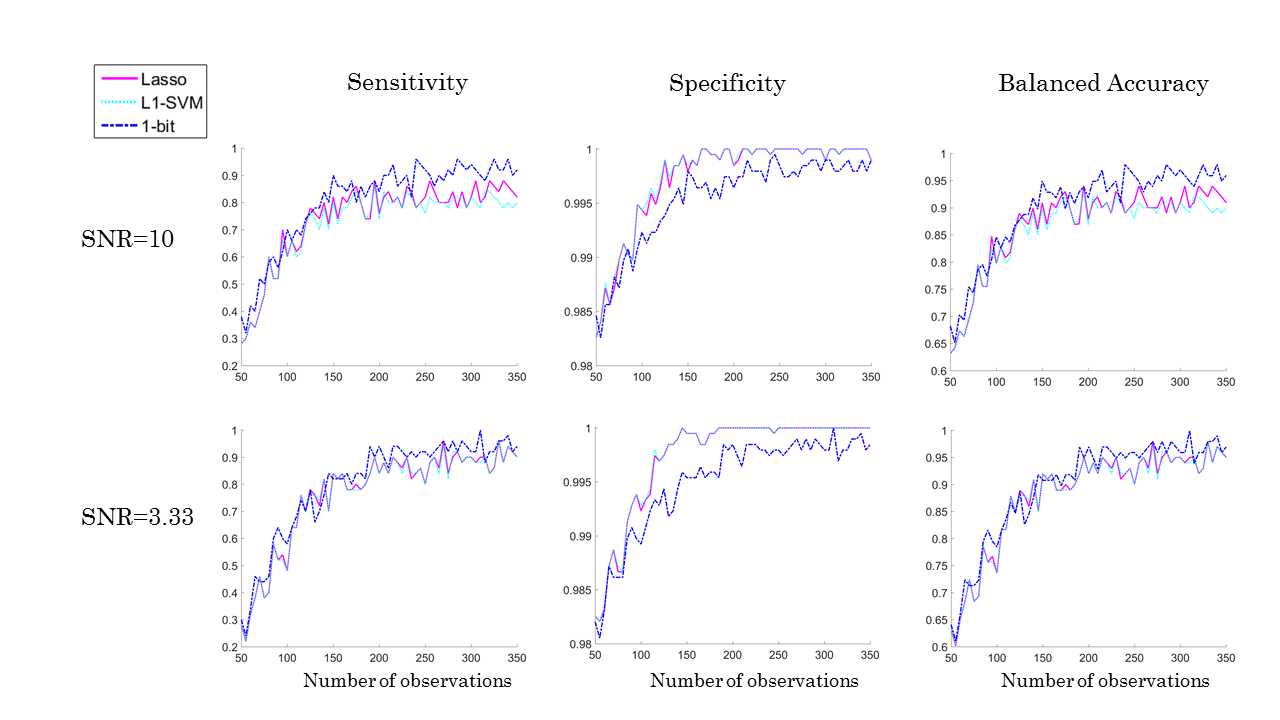}
		\caption{Comparison of numerical results for SPA ($= \text{$1$-bit CS}$), Lasso, and $\ell_1$-SVM on the data-set DS2 with SNR = 10 and 3.33 showed in the respective row.}
		\label{fig::DS2}

\end{figure}

For the sake of statistical stability, each experiment was repeated $10$-times. The averaged results are presented in the Fig.~\ref{fig::DS1}. We can see that SPA ($= \text{$1$-bit CS}$) performs better than the $\ell_1$-SVM or Lasso with regard to the capability of recognizing the true positive features (sensitivity in Fig. \ref{fig::DS1}). In our setting, if one method fails to select 5 condition positive peaks because one of them was selected twice, and the other method selects exactly the same 4 peaks and one false positive in addition, the specificity penalizes only the latter one. But effectively, both cases are suboptimal, since only the 5 positive peaks together can predict the class correctly.
This effect is reflected by a smaller value of specificity of the 1-bit approach comparing to the specificity of other two methods for data-sets with less than 300 spectra (column 2 in Fig. \ref{fig::DS1}). However, this also implies that SPA performs sligtly worse in rejecting true negatives than the other two approaches. The average results for balanced accuracy are visualized in the third column of Fig.~\ref{fig::DS1}. We observe that SPA outperforms the other two methods and even achieves 100\% accuracy with relatively few observations. With further decreasing SNR the performance of the three algorithms becomes more similar. %for data sets with fewer observations than the other two algorithms.
Fig.~\ref{fig::DS2} shows the numerical outcomes for the data-set DS2.
The non-trivial correlation structure of DS2 eventually leads to a slight drop of sensitivity and accuracy for SPA (compared to DS1), whereas the performance of the other two methods essentially remains unaffected. As before with further decreasing SNR the performance of the three algorithms becomes more similar in terms of sensitivity and balanced accuracy.

%%%%%%%%%%%%%%%%%%%%%%%%%%%%%%%%%%%%%%%%%%%%%%%%%%%%%%%%%%%%%%%%%%%%%%%%%%%%%%

%%%%%%%%%%%%%%%%%%%%%%%%%%%%%%%%%%%%%%%%%%%%%%%%%%%%%%%%%%%%%%%%%%%%%%%%%%%%%%
% \clearpage
\section{Experimental Results: Analyzing Real-World MALDI-TOF MS Data}
\label{section::experiments_real}
\begin{figure}
	\includegraphics[scale = 0.40]{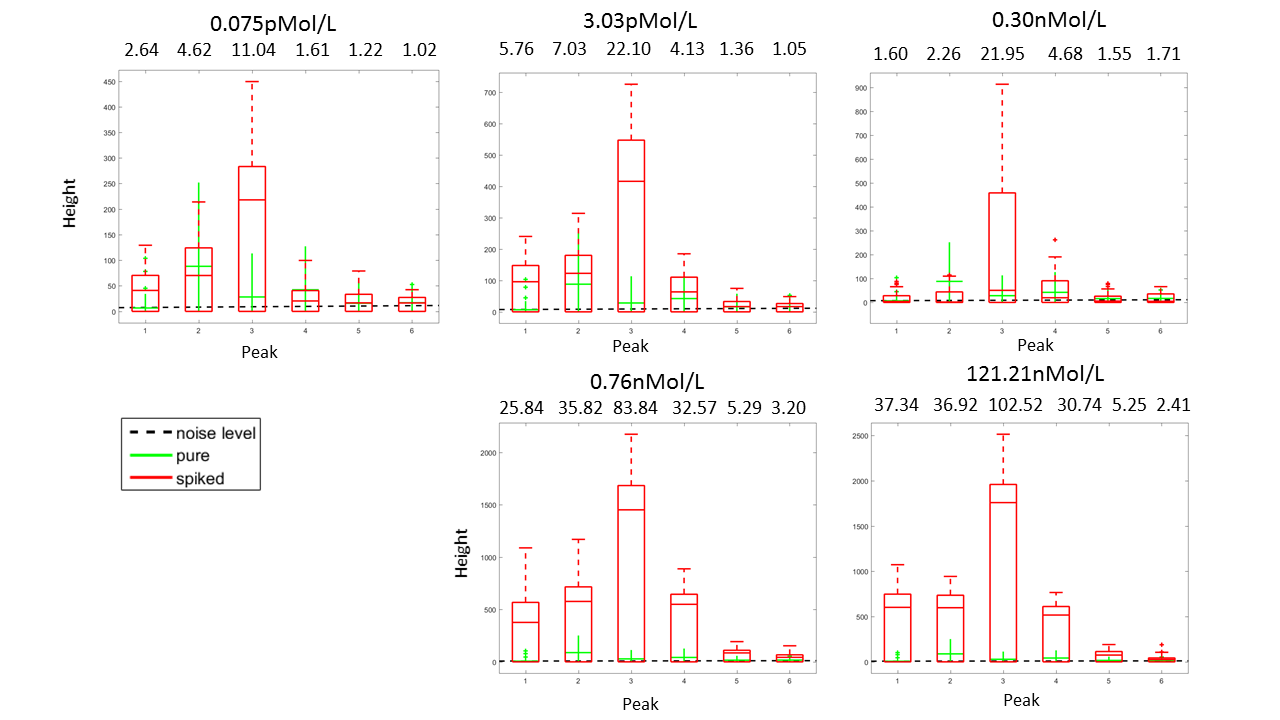}
	\caption{The height of true signals (6 spiked in peaks) comparing to the height of noise and height of the corresponding values in the pure data-set. Signal-to-noise ratio, which was calculated as the ratio of median of spiked-in signals and the estimated level of noise is shown above the corresponding peaks.}
	\label{fig::SNR}
	%	\captionsetup[subfigure]{justification=centering}	
	%	\centering
	%		\begin{subfigure}[b]{0.3\textwidth}
	%			\includegraphics[scale = 0.23]{figures/experiments/noise_vs_peak_intensity1.png}
	%			\subcaption{}
	%		\end{subfigure}	
	%		\begin{subfigure}[b]{0.3\textwidth}
	%			\includegraphics[scale=0.23]{figures/experiments/noise_vs_peak_intensity2.png}
	%			\subcaption{}
	%			\label{}
	%		\end{subfigure}	
	%		\begin{subfigure}[b]{0.3\textwidth}
	%			\includegraphics[scale=0.23]{figures/experiments/noise_vs_peak_intensity3.png}
	%			\subcaption{}
	%			\label{}
	%		\end{subfigure}	
	%		\begin{subfigure}[b]{0.3\textwidth}
	%			\includegraphics[scale=0.23]{figures/experiments/noise_vs_peak_intensity4.png}
	%			\subcaption{}
	%			\label{}
	%		\end{subfigure}	
	%		\begin{subfigure}[b]{0.3\textwidth}
	%			\includegraphics[scale=0.23]{figures/experiments/noise_vs_peak_intensity5.png}
	%			\subcaption{}
	%			\label{}
	%		\end{subfigure}	
	%		\caption{text}
	
\end{figure}

In this section, we present results of SPA, Lasso, and $\ell_1$-SVM for analyzing real-world mass-spectrometry data and compare them to the MALDIquant proteomics analysis workflow \cite{MALDIquant}. All data was acquired in our earlier studies \cite{Kratzsch2005, Fiedler2009}. It was approved by the local ethics committees and fulfils the requirements of the Helsinki declaration. All subjects gave informed consent to participate in the study. We will demonstrate the performance of our method on two data-sets:

\begin{itemize}
\item
	\emph{Spiked Data}: 
The spiked data-set is a labelled ground-truth data-set containing \emph{control} (e.g. healthy) and \emph{case} (e.g. diseased) mass spectra where the true labels are known. It is created from human blood samples\footnote{Blood serum of $16$ apparently healthy individuals from a clinical study (\cite{Kratzsch2005}) was used.} which were either unchanged (control group) or in which a protein-mix has been mixed (spiked) into (case group). In order to simulate different strength of an effect caused e.g. by a disease, we further sub-divided the \emph{case} group into five sub-groups where the amount of spiked-in proteins is increasing.
The five volumes in the case sub-groups were spiked with the following concentrations of the protein mix\footnote{Protein calibration standard mix Part No.: 206355 \& 206196) from Bruker Daltronics (Leipzig, Germany)}: 0.075pMol/L, 3.03pMol/L, 0.30nMol/L, 0.76nMol/L and 121.21nMol/L. This mix contains the hormones Angiotensin, ACTH, clip 18-39, Substance P and the cell protein Ubiquitin. The peptide mix was added before sample pre-treatment and 64 spectra were measured due to 4-fold spotting (technical replicates). Mass spectra were acquired using the protocol described in the supplementary material (S1). Each volume corresponds to a data-set. What differentiates the data-sets are the amplitudes of the 6 spikes resulting from the added substances. The signal-to-noise ratio of the spiked-in peaks is shown in the Fig. \ref{fig::SNR}\footnote{The power of noise for each of the 5 analyzed data-sets is estimated as an average of intensity of noise of the observations using median absolute deviation.}.
\item
	\emph{Pancreas Cancer Data (P. CA)}: A total of 120 patients with pancreatic cancer and controls were recruited for this study \cite{Fiedler2009}. For the discovery study sera were obtained from two different clinical centres (University Hospital Leipzig (UHL, set L) and Heidelberg (UHH, set H)) as described in the supplementary material (S1). Note that each acquired spectrum has been assigned a class-label, i.e., healthy or diseased. So, the health status of the training samples is known in advance (supervised learning).
\end{itemize}
% \vspace{0.3cm}

Baseline removal was performed on the raw MS data using TopHat filtering (\cite{Sauve2004}). In particular, no additional calibration or noise reduction steps have been applied. More information on the data and sample preparation can be found in the supplementary material (S1).

\subsubsection*{The \emph{missing-data} problem}
When dealing with data coming from measurements of, say, a Mass Spectrometer instrument, the so called \emph{missing-data problem} usually occurs. This means that the instrument failed to give measurements for some of the measured masses, usually due to the stochastic nature of the process happening inside the device. Due to the smoothing step in our algorithm and the arguments of e.g. Rubin et al. (\cite{Ru1976}) this problem can be mainly ignored in our case for identifying the relevant features. However, this does not necessarily hold for the classification step, i.e. applying the identified sparse classifier to an unknown data-set. In this scenario, where data is missing in an unknown sample, there are basically two options: (1) applying a method for inferring the missing data or (2) stopping the classification and return an error message to the user. In this work we decided to follow the latter approach, since inferring missing data is not in the scope of this paper\footnote{The interested reader might find a good starting point about this topic in these two reviews \cite{Pi2001,Sch1998}} but is an unarguable crucial point in any data analysis pipeline and should depend on the actual use-case.

% \vspace{0.3cm}

\subsubsection*{Accuracy vs. Number of Features}

We performed these experiments with respect to the same evaluation categories as in the case of simulated data. Note that the normalization and standardization as described in Section \ref{section::AlgorithmicDetails} were applied as preprocessing steps in each of the methods. Similarly, a hard thresholding as described in \eqref{eq:hardthres} was applied to all classifiers estimated by the examined algorithms.

%The spiked data-sets were declared as condition positives while the pure - unspiked spectra were declared condition negatives. This procedure provides us with the binary labels that are a prerequisite for applying SPA. 
For the each of the algorithms, we are testing the performance of the obtained classifiers learned on the pure data-set which corresponds to the condition negative class and one spiked data-set at a time corresponding to the condition positive class.
% Similar to the simulation study, we know the ground-truth features. The 6 spiked-in peaks (see the data-set description above) are precisely the discriminating features that we are interested in and their location is the best information that we have a priory about the location of biomarkers. The main question in our experiments is therefore, how successful each of the algorithms is in detecting the 6 peaks that were initially spiked. Asking for the detection of the correct spikes is more natural here and is an even more difficult task than classifying.

%\begin{table}[]
%	\centering	
%	\begin{tabular}{|c|c|c|c|c|c|}
%		\hline
%		data-set      & 0.075pMol/L & 3.03pMol/L & 0.30nMol/L & 0.76nMol/L & 121.21nMol/L \\ \hline
%		SNR of peak 1 & 2.64        & 5.76       & 1.60       & 25.85      & 37.34        \\ \hline
%		SNR of peak 2 & 4.64        & 7.03       & 2.26       & 35.82      & 36.92        \\ \hline
%		SNR of peak 3 & 11.04       & 22.10      & 21.95      & 83.84      & 102.52       \\ \hline
%		SNR of peak 4 & 1.61        & 4.13       & 4.68       & 32.57      & 30.74        \\ \hline
%		SNR of peak 5 & 1.22        & 1.36       & 1.55       & 5.29       & 5.25         \\ \hline
%		SNR of peak 6 & 1.02        & 1.05       & 1.71       & 3.20       & 2.41         \\ \hline
%	\end{tabular}
%	\caption{Signal-to-noise ratio of spiked-in peaks.}
%	\label{tab::SNR}
%\end{table}

The results of the classifier with $6$ non-zeros on the spiked data-set are shown in Table~\ref{tab:main_results_spiked}. The main question in these experiments is how successful each of the algorithms is in detecting the $6$ peaks that were initially spiked (see the data-set description above). We can see that the values of sensitivity for SPA are at least as high as those of the other methods, which implies that the approach of $1$-bit CS is very competitive in this situation and mostly achieves the best detection rate. However, the relatively poor performance of all the algorithms on the spiked data-set can be explained by the nature of the data. Since the peptide mix was added to the blood samples before acquiring the mass spectra, the spiked peaks are not always present in all the resulting mass spectra in the positions where we expect to find them. There exist data-sets for which all the mass spectra failed to exhibit certain spiked peaks at their expected locations as can be seen in the Fig. \ref{fig::SNR}. Thus, we cannot expect any of the algorithms to find these missing peaks. Nonetheless, there is still a chance to build a reliable fingerprint out of the remaining spikes while there is no chance to detect the missing spikes because the data-set is not rich enough to represent it. On the other hand, this spiked data-set combines the advantages of both simulated and clinical data, since the positions of the desired biomarkers are known in advance while their representative behavior in the spectra is quite realistic.

\begin{figure}[!ht]
	
	\begin{subfigure}[b]{0.45\textwidth}
		\includegraphics[scale = 0.40]{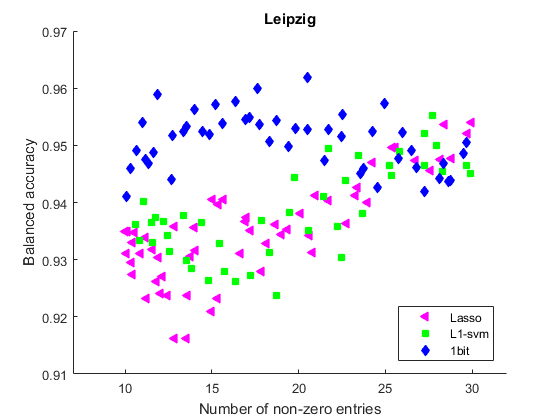}
		\label{subfig::Leipzig}
	\end{subfigure}
	\begin{subfigure}[b]{0.45\textwidth}
		\includegraphics[scale = 0.40]{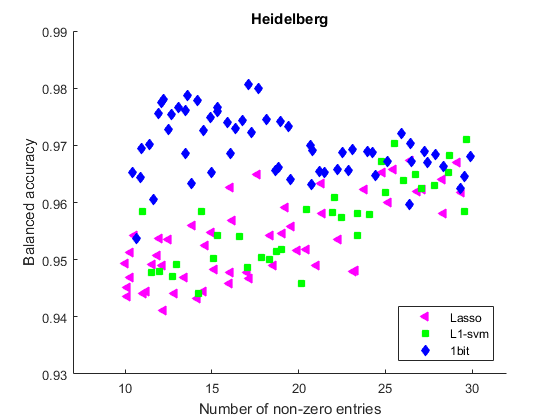}
		\label{subfig::Heidelberg}
	\end{subfigure}
	
	\caption{Accuracies of sparse classifiers from SPA, Lasso, and $\ell_1$-SVM on the real pancreatic cancer data-sets. While Lasso and $\ell_1$-SVM achieve better classification accuracy with increasing number of features, SPA is particularly well suited for the ``very-sparse regime'' where only few features ($<20$) are used for classification.}
	\label{fig::real}
	
\end{figure}

In contrast to that in the case of pancreas cancer data-sets, we do not know the true-positive feature positions. Consequently, we can only rely on the classification performance of the obtained sparse classifiers by each of the algorithms. To evaluate the reliability of our results, for each of the methods, we have employed the cross-validation scheme as described in the Algorithm \ref{algo:CV} with the number of folds $K$ set to 5. In order to ensure statistical stability, each experiment was repeated 10-times. Fig.\ref{fig::real} shows the average results over 10 repetitions.

		\algobox[Cross-Validation of Classification Performance]{.95\textwidth}{\label{algo:CV}
			
			\emph{Input:} Raw data $(x_1,y_1),...,(x_n,y_n)\in \R^d\times \{-1,+1\};$
			Number of CV-folds $K$;
			
			\emph{Output:} Classification accuracy $Acc\in [0,\ 1]$;
			Average sparsity $\#F$ (number of selected features)

			\begin{enumerate}
				\item[1:]
				Split the sample set $\{1,...,n\}$ randomly into $K$ disjoint folds $P_1, ..., P_K$ of (almost) equal size. \vspace{6pt}\\ 
				\emph{For each fold $k\in{1,...,K}$ perform the following steps 2-4:}
				\item[2:]
				Compute the feature vector $\omega^k$ employing the desired method on the samples of $\bigcup_{ k'\in \{1,...,K\} \setminus \{k\} }P_{k'}$
				\item[3:]
				Dimension reduction as described in Sec. \ref{subsec:dim-reduction} (step 6). Project all spectra onto $\operatorname{supp}(\omega^k)$ and put $\#F^k := ||\omega^k||_0.$
				\item[4:]
				Classification of $P_k$: Use the projected samples of $\bigcup_{ k'\in \{1,...,K\} \setminus \{k'\} }P_{k'}$ to predict the labels of the spectra $P_k$ by an ordinary SVM\footnotemark. Denote the prediction accuracy by $Acc^k$.
				\item[5:]
				Compute the average accuracy $\operatorname{Acc}:= \frac{1}{K}\sum_{k=1}^{K}Acc^k$ and average sparsity $\#F:= \frac{1}{K}\sum_{k=1}^{K} \#F^k.$
			\end{enumerate}
		}
		\footnotetext{Here, the standard MATLAB implementation of SVM was used.}

% Instead, for each of the methods, after computing the sparse classifier on one part of the data-set (using a 5-fold cross-validation schema) and dimensionality reduction by projecting data onto the selected features, classification of the remaining part of the data-set using linear SVM was performed\footnote{Here, the standard MATLAB implementation of SVM was used.} using the provided class-labels of the data.
In order to ensure statistical stability, each experiment was repeated $10$-times. Fig.~\ref{fig::real} shows the results.
Note that our results show that accurate predictions are already possible with a very few features, so that the assumption of small disease fingerprint seems to hold for this data-set. Furthermore, it can be seen that SPA is especially well suited for situations where a sparse classifier (containing only very few features) is preferred. This is appealing because fewer features enable an easier interpretation of the actual components of a potential disease fingerprint. Moreover, follow-up experiments that often involve an individual treatment of each component (e.g., potential biomarkers) would become much less costly.
Note that in the non-sparse region with more than $30$ features selected, it is not meaningful to relate the achieved accuracy to the quality of the learned feature vector due to the small sample size. The considered algorithms assume the underlying fingerprint to be sparse. This assumption usually does not fully hold in practice. Therefore, we cannot expect that a learned feature vector achieves perfect classification. The classification accuracy should be therefore considered as an indicator of how well our model assumption of the sparse fingerprint fits to the unknown ground-truth. If we let the algorithms operate out of the region for which they have been designed for, we may achieve indeed a higher accuracy, but this is probably a consequence of overfitting. And even more importantly, the learned feature vector (model) is not reliable anymore.

%%%%%%%%%%%%%%%%%%%%%%%%%%%%%%%%%%%%%%%%%%%%%%%%%
\subsubsection*{Best Classifier}

Apart from that, we are interested in the performance of the best sparse classifier (i.e. small number of features) found by each of the algorithms (SPA, Lasso, $\ell_1$-SVM). For all learned classifiers with $10$ to $30$ non-zero components, Table~\ref{tab:main_results_real} presents those with the best classification accuracy. Furthermore, we also considered a typical analysis pipeline (MALDIQuant) to see how the ``purely-data-based'' approaches (SPA, Lasso, $\ell_1$-SVM) compare to a model-based approach\footnote{By ``model-based'' we mean that specific model assumptions on the data are made and exploited, such as noise-structure for denoising or Gaussian-shaped structures for peak detection.}. In Table~\ref{tab:main_results_real}, it can be seen that SPA provides the sparsest solutions while achieving competitive results with respect to sensitivity and specificity at the same time. Lasso and $\ell_1$-SVM select almost the same features and therefore perform similarly.
On the other hand, MALDIQuant selects the features based on a prior model-based peak detection followed by a feature selection based on shrinkage diagonal discriminant analysis (\cite{DDA}). But however, it still performs worst on the UHL data-set.

%%%%%%%%%%%%%%%%%%%%%%%%%%%%%%%%%%%%%%%%%%%%%%%%%%%%%%%%%%%%%%%%%%%%%%%%%%%%%%
%% Results of SPIKED data
%%%%%%%%%%%%%%%%%%%%%%%%%%%%%%%%%%%%%%%%%%%%%%%%%%%%%%%%%%%%%%%%%%%%%%%%%%%%%%
% \clearpage

\begin{landscape}

\newcolumntype{L}{>{\arraybackslash}m{2.5em}}

\begin{table}
	\begin{filecontents*}{results/results_for_paper_spiked-data.csv}
concentration;tp-spa;sens-spa;spec-spa;acc-spa;tp-lasso;sens-lasso;spec-lasso;acc-lasso;tp-l1svm;sens-l1svm;spec-l1svm;acc-l1svm
0.075pMol/L;2;0.333333333;0.999905605;0.666619469;1;0.166666667;0.999976401;0.583321534;1;0.166666667;1;0.583333333
3.03pMol/L;4;0.666666667;0.999952802;0.833309735;2;0.333333333;1;0.666666667;1;0.166666667;1;0.583333333
0.30nMol/L;2;0.333333333;0.999905605;0.666619469;1;0.166666667;1;0.583333333;1;0.166666667;1;0.583333333
0.76nMol/L;2;0.333333333;0.999905605;0.666619469;2;0.333333333;0.999952802;0.666643068;2;0.333333333;0.999952802;0.666643068
121.21nMol/L;3;0.5;0.999929204;0.749964602;2;0.333333333;0.999976401;0.666654867;2;0.333333333;0.999976401;0.666654867
\end{filecontents*}
\pgfplotstabletypeset[
	%	skip first n=1,
	font={\small},
	col sep=semicolon,
	string type,
	every head row/.style={%
		before row={\hline
			&  \multicolumn{8}{c}{SPA} & \multicolumn{8}{c}{Lasso} & \multicolumn{8}{c}{$\ell_1$-SVM}\\
		},
		after row=\hline
	},
	every last row/.style={after row=\hline},
	every odd row/.style={before row={\rowcolor[gray]{0.9}}},
	% dataset
	columns/concentration/.style={column name=Concentration, column type=l|},
	% spa
	columns/tp-spa/.style={column name=TP$^{[a]}$, column type=l, fixed, fixed zerofill, dec sep align, precision=0},
	columns/sens-spa/.style={column name=Sens$^{[b]}$, column type=l, fixed, fixed zerofill, dec sep align, precision=3},
	columns/spec-spa/.style={column name=Spec$^{[c]}$, column type=l, fixed, fixed zerofill, dec sep align, precision=3},
	columns/acc-spa/.style={column name=B. Acc$^{[d]}$, fixed, fixed zerofill, dec sep align={l|}, precision=3, column type/.add={}{|}},
	% Lasso
	columns/tp-lasso/.style={column name=TP, column type=l, fixed, fixed zerofill, dec sep align, precision=0},
	columns/sens-lasso/.style={column name=Sens, column type=l, fixed, fixed zerofill, dec sep align, precision=3},
	columns/spec-lasso/.style={column name=Spec, column type=l, fixed, fixed zerofill, dec sep align, precision=3},
	columns/acc-lasso/.style={column name=B. Acc, fixed, fixed zerofill, dec sep align={l|}, precision=3, column type/.add={}{|}},
	% l1-svm
	columns/tp-l1svm/.style={column name=TP, column type=l, fixed, fixed zerofill, dec sep align, precision=0},
	columns/sens-l1svm/.style={column name=Sens, column type=l, fixed, fixed zerofill, dec sep align, precision=3},
	columns/spec-l1svm/.style={column name=Spec, column type=l, fixed, fixed zerofill, dec sep align, precision=3},
	columns/acc-l1svm/.style={column name=B. Acc, fixed, fixed zerofill, dec sep align, precision=3},
	] {results/results_for_paper_spiked-data.csv}
	\caption{This table shows the main results comparing the feature selection benchmarks of our approach with Lasso and $\ell_1$-SVM on the spiked data-set. Given results correspond to the average results over 10 repetitions of the classifier with 6 non-zero values.\\~\\}
		$[a]$ \textbf{TP}: Number of spiked peaks that are correctly detected;
		$[b]$ \textbf{Sens}: Sensitivity in detecting spiked peaks ($TP/(TP + FN)$);
		$[c]$ \textbf{Spec}: Specificity in detecting spiked peaks ($TN/(FP + TN)$);
		$[d]$ \textbf{B. Acc}: Balanced Accuracy ($ \frac{\operatorname{sens.} + \operatorname{spec.}}{2}$);
		\label{tab:main_results_spiked}
\end{table}

%%%%%%%%%%%%%%%%%%%%%%%%%%%%%%%%%%%%%%%%%%%%%%%%%%%%%%%%%%%%%%%%%%%%%%%%%%%%%%
%% Results of REAL data
%%%%%%%%%%%%%%%%%%%%%%%%%%%%%%%%%%%%%%%%%%%%%%%%%%%%%%%%%%%%%%%%%%%%%%%%%%%%%%

\begin{table}
\begin{filecontents*}{results/results_for_paper_real-data.csv}
dataset;feature-count-spa;sens-spa;spec-spa;acc-spa;feature-count-lasso;sens-lasso;spec-lasso;acc-lasso;feature-count-l1svm;sens-l1svm;spec-l1svm;acc-l1svm;feature-count-mquant;sens-mquant;spec-mquant;acc-mquant
P. CA - UHL;20.48;0.975;0.948833;0.961917;29.94;0.96875;0.939333;0.954042;27.72;0.96375;0.946583;0.955167;21;0.888;0.888;0.888
P. CA - UHH;17.1;0.986286;0.975;0.980643;26.46;0.966095;0.96875;0.967423;29.68;0.965619;0.97625;0.970935;17;0.975;0.963;0.969
\end{filecontents*}
\pgfplotstabletypeset[
%	skip first n=1,
	font={\small},
    col sep=semicolon,
    string type,
    every head row/.style={%
        before row={\hline
          &  \multicolumn{7}{c}{SPA} & \multicolumn{7}{c}{Lasso} & \multicolumn{7}{c}{$\ell_1$-SVM} & \multicolumn{7}{c}{Maldi-Quant}\\
        },
        after row=\hline
    },
    every last row/.style={after row=\hline},
   	every odd row/.style={before row={\rowcolor[gray]{0.9}}},
    % dataset
    columns/dataset/.style={column name=Dataset, column type=l|},
    % spa
    columns/feature-count-spa/.style={column name=Feat.$^{[e]}$, column type=L},
    columns/sens-spa/.style={column name=Sens$^{[f]}$, column type=l, fixed, fixed zerofill, dec sep align, precision=3},
    columns/spec-spa/.style={column name=Spec$^{[g]}$, column type=l, fixed, fixed zerofill, dec sep align, precision=3},
    columns/acc-spa/.style={column name=B. Acc$^{[h]}$, fixed, fixed zerofill, dec sep align={l|}, precision=3, column type/.add={}{|}},
    % lasso
    columns/feature-count-lasso/.style={column name=Feat., column type=L},
    columns/sens-lasso/.style={column name=Sens, column type=l, fixed, fixed zerofill, dec sep align, precision=3},
    columns/spec-lasso/.style={column name=Spec, column type=l, fixed, fixed zerofill, dec sep align, precision=3},
    columns/acc-lasso/.style={column name=B. Acc, fixed, fixed zerofill, dec sep align={l|}, precision=3, column type/.add={}{|}},
    % l1-svm
    columns/feature-count-l1svm/.style={column name=Feat., column type=L},
    columns/sens-l1svm/.style={column name=Sens, column type=l, fixed, fixed zerofill, dec sep align, precision=3},
    columns/spec-l1svm/.style={column name=Spec, column type=l, fixed, fixed zerofill, dec sep align, precision=3},
    columns/acc-l1svm/.style={column name=B. Acc, fixed, fixed zerofill, dec sep align={l|}, precision=3, column type/.add={}{|}},
    % maldi-quant
    columns/feature-count-mquant/.style={column name=Feat., column type=L},
    columns/sens-mquant/.style={column name=Sens, column type=l, fixed, fixed zerofill, dec sep align, precision=3},
    columns/spec-mquant/.style={column name=Spec, column type=l, fixed, fixed zerofill, dec sep align, precision=3},
    columns/acc-mquant/.style={column name=B. Acc, column type=l, fixed, fixed zerofill, dec sep align, precision=3},
    ] {results/results_for_paper_real-data.csv}
\caption{This table shows the main results comparing the feature selection benchmarks of our approach with Lasso, $\ell_1$-SVM, and Maldi-Quant. These are averages over $10$ repetitions of a $5$-fold cross-validation.  Note that these results have been calculated based on the highest accuracy criterion for all classifiers with between $10$ and $30$ selected features. This particularly means that better accuracy values might be achieved for the individual methods if less sparse feature vectors would be allowed. For more details see text. \\ ~\\
$[e]$ \textbf{Feat.}: Number of features;
$[f]$ \textbf{Sens}: Sensitivity ($TP/(TP + FN)$);
$[g]$ \textbf{Spec}: Specificity ($TN/(FP + TN)$);
$[h]$ \textbf{B. Acc}: Balanced Accuracy ($ \frac{\operatorname{sens.} + \operatorname{spec.}}{2}$);
}
\label{tab:main_results_real}
\end{table}

\end{landscape}

%%%%%%%%%%%%%%%%%%%%%%%%%%%%%%%%%%%%%%%%%%%%%%%%%
\subsubsection*{Medical Interpretation of Results}

Pancreatic cancer is not only a common and increasingly frequent \cite{Yeo2012}, but also
still a fatal disease, with a survival rate of 3-5\% five years after diagnosis
\cite{Michl2006}. The conventional tumor marker, Carbohydrate Antigen 19-9 (CA19-9), as
a blood group antigen not present in a significant proportion of the patients
\cite{Leichtle2015}, shows insufficient diagnostic sensitivity and specificity (AUC 0.71),
even in combination with the second-line tumor marker Carcinoembryonic Antigen
(CEA, combined AUC 0.75) \cite{Reitz2015}. The need for better markers for screening and
differential diagnosis is evident, as panceratic carcinoma would be principally
curable if detected and identified very early in the course of the disease.
Along with the emerging ``-omics''-technologies great hope was rised to find
tumor-specific peptides or metabolic alterations to increase sensitivity and
specificity of early and differential diagnostics, and several combinatory
marker models could be identified by proteomics \cite{Fiedler2009} and metabolomics \cite{LeiCegCon2013}. Pancreatic carcinoma is a complex disease - it affects the metabolism as a
whole (e.g. the so-called Warburg effect) \cite{Zhou2012b}, but also alters proteolytic
activity \cite{Brand2011}. Therefore, it might be na\"{\i}ve to expect a single marker capable
to indicate presence, progression and exact type of the malignancy at once \cite{Leichtle2013}
9– it might even be overly reductionistic to attribute these capabilities to a
single model, even if it consists of several entities measured by different ``-
omics'' technologies \cite{Leichtle2015}. As Raftery states ``basing inferences on a single
``best'' model as if the single selected model were true ignores model uncertainty,
which can result in underestimating uncertainty about quantities of interest''
\cite{Raftery1997}, and the larger the ``-omics'' data-sets grow, the larger is the
`probability, that there is not one ``single best'' predictive marker model, but
instead several with comparable selectivity \cite{Leichtle2013}. And it is very reasonable
to assume that, even on the same data-set, different algorithms might favor
different models consisting of different feature sets and bring forth completely
different results, when only the best differentiating models are regarded. For
an in-depth comparison of the validity of the results of different algorithms,
the underlying peak features should also be taken into account, and similarities
in the selected features corroborate the algorithms superimposed on them. In the
case of our study, we have the great advantage, that the same data-set was
evaluated in three different studies: the principal one by Fiedler et al. \cite{Fiedler2009}, a subsequent BinDA-algorithm-based manuscript by Gibb and Strimmer published
recently \cite{Gibb2015}, and the present one. Fiedler et al. \cite{Fiedler2009} identified one
discriminating peptide, Platelet Factor 4 (m/z \textit{3884}, identified in italics,
double hits in bold) within four discriminating peaks (m/z 3194, \textit{3884}, 4055,
and \textbf{5959}). The 30 most differential peaks in Gibb et al. \cite{Gibb2015} were m/z
4495, 8868, 8989, 1855, 4468, 8937, 2023, 1866, 5864, 5946, 1780, 2093, \textbf{5906},
\textbf{5960}, 8131, 1207, 4236, 2953, 9181, 1021, \textbf{\textit{1466}}, 4092, 4251, 5005, 8184, 1897, 3264, 2756, 6051, and 1264, with m/z 8937 identified as pancreatic progenitor
cell differentiation and proliferation factor-like protein. m/z 3884 could not
be identified as discriminating marker (while it might play a role in pancreatic
carcinoma nonetheless \cite{Poruk2010}), whereas m/z \textbf{\textit{1466}} can be attributed to a
fragment of fibrinopeptide A (DSGEGDFLAEGGGVR), as previously described in tumor
samples \cite{Villanueva2006}. In the present study, the peaks  m/z \textbf{\textit{1464}}, 1546, 1944, \textbf{5904}, 1619, 4209, and 2662 could be identified as discriminating features. The
slight mass shift of about 2 Da for m/z \textbf{\textit{1464}} / \textbf{\textit{1466}} and \textbf{5904} / \textbf{5906} is
probably arising from different peak preprocessing procedures, peaks are wide
enough to tolerate this deviation. Further investigations and the application of
further algorithms on the same data-set are highly likely to yield a similar,
partially overlapping set of features, each with a comparable discriminating
power (Fiedler et al. \cite{Ceglarek2009} AUC$_{[3884/(CA19-9*CEA)]}$ 1.0;
Gibb et al. \cite{Gibb2015} in a 5-feature model: accuracy of 0.96,
sensitivity of 0.96, specificity of
0.97, positive predictive value of 0.97 and negative predictive value of 0.95;
the present study accuracy$_{[UHL]}$ 0.96, sensitivity$_{[UHL]}$ 0.97, specificity$_{[UHL]}$ 0.95 and accuracy$_{[UHH]}$ 0.98, sensitivity$_{[UHH]}$ 0.99, specificity$_{[UHH]}$ 0.97. This also corresponds to a
recently published comparable study investigating a glycoprotein marker panel
(AUC 0.95) \cite{Nie2014}. Biomarkers for clinical diagnostics comprise a wide field
of applications (e.g. population-wide screening, early diagnostics,
characterization, treatment guidance, efficacy and toxicity monitoring,
prognosis, susceptibility estimation and many more) \cite{Leichtle2015}, each with special
requirements for sensitivity and specificity, that are only partially condensed
in the AUC as an overall selectivity measure \cite{Leichtle2013}. Especially for screening
purposes, sensitivity is extremely important \cite{LeiCegCon2013}, and clinically applied
tests e.g. for newborn screening frequently surpass the 0.99 hallmark \cite{Ceglarek2009}.
Compared with the conventional, ``not-for-screening'' marker CA19-9, the SPA-based
model shows considerable improvement, however there is still a big gap to
screening suitability, which in the next years might be bridged by improved
sensitivity of new instrumentation, refined algorithms (as the SPA), and combination with
other ``markers'' from the ``big data'' field, enabling a more
holistic view – not only of the disease, but also of the affected patient\cite{Leichtle2015}.

%%%%%%%%%%%%%%%%%%%%%%%%%%%%%%%%%%%%%%%%%%%%%%%%%%%%%%%%%%%%%%%%%%%%%%%%%%%%%%

% \pagestyle{empty}
% \begin{landscape}
%
% \begin{figure}
%	\includegraphics[scale = 0.31]{figures/experiments/maldiquant1.png}
%	\caption{The MALDIquant workflow: the raw spectrum (top), spectrum after preprocessing (middle), detected peaks  (bottom).}
% \end{figure}

% \end{landscape}
% \pagestyle{plain}

%%%%%%%%%%%%%%%%%%%%%%%%%%%%%%%%%%%%%%%%%%%%%%%%%%%%%%%%%%%%%%%%%%%%%%%%%%%%%%
%\clearpage
%\section{\methodname \hspace{0.1cm} for Clustering and Classification}
%\label{section::classification}
%\input{sparseProteomicsAnalysis_section_clustering-classification}
%%%%%%%%%%%%%%%%%%%%%%%%%%%%%%%%%%%%%%%%%%%%%%%%%%%%%%%%%%%%%%%%%%%%%%%%%%%%%%

%%%%%%%%%%%%%%%%%%%%%%%%%%%%%%%%%%%%%%%%%%%%%%%%%%%%%%%%%%%%%%%%%%%%%%%%%%%%%%
% \clearpage
\section{Conclusion}
%%%%%%%%%%%%%%%%%%%%%%%%%%%%%%%%%%%%%%%%%%%%%%%%%%%%%%%%%%%%%%%%%%%%%%%%%%%%%%
In this work, we introduced Sparse Proteomics Analysis (SPA), a new framework for the analysis of proteomics data, generated by mass spectrometry measurements. The framework solves the problem of selecting a minimum set of features from high-dimensional data in a situation where only relatively few measurements are available. Our approach particularly allows for a bio-medical interpretation and enables a classification of unknown samples. This is done by formulating and solving a regularized optimization problem, using ideas from $1$-bit compressed sensing combined with several generic pre- and postprocessing steps. We have shown by several numerical experiments that SPA can indeed compete with standard (and widely used) algorithms as well as with a specifically adapted method for the analysis of proteomics data (MALDIQuant). In the ``very-sparse'' situation, it has even turned out that SPA outperforms the other approaches with respect to prediction accuracy.

%%%%%%%%%%%%%%%%%%%%%%%%%%%%%%%%%%%%%%%%%%%%%%%%%%%%%%%%%%%%%%%%%%%%%%%%%%%%%%
% \clearpage
\section*{Competing Interests}
The authors declare that they have no competing interests.

\section*{Author's Contributions}
Conceived and designed the experiments: TC, MG, JV, GK, and CS. Performed the experiments: NC, TC, and NW. All authors contributed to writing the paper. All authors read and approved the final manuscript.

\section*{Acknowledgments}
JV was supported by the ERC CZ grant LL1203 of the Czech Ministry of Education. TC, MG, NC, JV, GK and CS were supported by the Einstein Center for Mathematics Berlin (ECMath), project grant CH2, and by the DFG Research Center Matheon \emph{Mathematics for key technologies}, Berlin. TC and CS are supported by the German Ministry of Research and Education (BMBF) project Grant 3FO18501 (Forschungscampus MODAL). GK acknowledges support by the Einstein Foundation Berlin, by the Deutsche Forschungsgemeinschaft (DFG), and by the DFG Collaborative Research Center TRR 109 \emph{Discretization in Geometry and Dynamics}.

The authors are thankful to Irena Bojarovska for fruitful discussions and help conducting the experiments.

%%%%%%%%%%%%%%%%%%%%%%%%%%%%%%%%%%%%%%%%%%%%%%%%%%%%%%%%%%%%%%%%%%%%%%%%%%%%%%%%
%%%%%%%%%%%%%%%%%%%%%%%%%%%%%%%%%%%%%%%%%%%%%%%%%%%%%%%%%%%%%%%%%%%%%%%%%%%%%%%%

\clearpage
\section*{Supporting Information}

% Include only the SI item label in the subsection heading. Use the \nameref{label} command to cite SI items in the text.

%%%%%%%%%%%%%%%%%%%%%%%%%%%%%%%%%%%%%%%%%%%%%%%%%%%%%%%%%%%%%%%%%%%%%%%%%%%%%%
%%%%%%%%%%%%%%%%%%%%%%%%%%%%%%%%%%%%%%%%%%%%%%%%%%%%%%%%%%%%%%%%%%%%%%%%%%%%%%
\subsection*{S1 - Mass Spectrometry Data Generation}
\label{S1_Text}
\label{supp_material::data}

\subsubsection*{Chemicals, Standards and Consumables}
Gradient grade acetonitrile, ethanol, and HPLC-water were obtained from J.T. Baker (Phillipsburg, NJ, USA); p.a. trifluoroacetic acid (TFA) and acetone were purchased from Sigma-Aldrich (Taufkirchen, Germany). The peptide- and protein MALDI-TOF calibration standards I and $\alpha$-cyano-4-hydroxycinnamic acid (HCCA) were purchased from Bruker Daltonics (Bremen, Germany).
Automated magnetic bead preparations were performed using 96 well plates, TubePlates from Biozym (Hessisch Oldendorf, Germany), polypropylene tubes (low profile) from Abgene (Surrey, UK), and modular reservoir quarter modules from Beckman (Fullerton, USA). For sample storage 450 $\mu$L CryoTubesTM were purchased from Sarstedt (Nmbrecht, Germany). Multifly needle sets and polypropylene serum monovettes with clotting activators were also obtained from Sarstedt.

\subsubsection*{Peptidome Separation}
All serum samples of the discovery set were processed at one time and analyzed simultaneously to avoid procedure-dependent errors. The external validation set was prepared, processed and analyzed separately. Peptidome separation of the samples was performed using the ClinPro Tools profiling purification kits from Bruker Daltonics. Magnetic particles with defined surface functionalities (magnetic beadimmobilized metal ion affinity chromatography (MB-IMAC Cu), magnetic bead-hydrophobic interaction (MB-HIC C8) and weak cation exchange (MB-WCX)) were processed by the ClinPro Tools liquid handling robot according to the manufacturers protocol (Bruker Daltonics).
Serum specimens were thawed on ice for 30 min and immediately processed according to our
standardized protocol for serum peptidomics \cite{Baumann2005}.

\subsubsection*{Mass Spectrometry}
A linear MALDI-TOF mass spectrometer (Autoflex I, Bruker Daltonics) was used for the peptidome profiling. Daily mass calibration was performed using the standard calibration mixture of peptides and proteins in a mass range of 1-10 kDa. Mass spectra were recorded and processed using AutoXecute tool of the flexControl acquisition software (Ver. 2.0; Bruker Daltonics).
The settings were applied as follows: Ion source 1: 20 kV; ion source 2, 18.50 kV; lens, 9.00 kV; pulsed ion extraction, 120 ns; nitrogen-pressure, 2500 mbar. Ionization was achieved by a nitrogen laser ($\lambda$=337 nm) operating at 50 Hz. For matrix suppression a high gating factor with signal suppression up to 500 Da was used. Mass spectra were detected in linear positive mode.

%%%%%%%%%%%%%%%%%%%%%%%%%%%%%%%%%%%%%%%%%%%%%%%%%%%%%%%%%%%%%%%%%%%%%%%%%%%%%%
\subsubsection*{Baseline Removal}
\label{supp_material::preprocessing}
%%%%%%%%%%%%%%%%%%%%%%%%%%%%%%%%%%%%%%%%%%%%%%%%%%%%%%%%%%%%%%%%%%%%%%%%%%%%%%
The baseline is an exponential like offset dependent on the $m/z$ value (mass-to-charge; x-value). A baseline correction is performed to remove this rather low-frequency noise from the spectrum. We use a morphological TopHat filter to eliminate certain spatial structures within the signal, in our case the baseline. Note that this technique does not produce negative intensity values.

\nolinenumbers

%\section*{References}
% Either type in your references using
% \begin{thebibliography}{}
% Text
% \end{thebibliography}
%
% OR
%
% Compile your BiBTeX database using our plos2015.bst
% style file and paste the contents of your .bbl file
% here.
%

\bibliographystyle{plos2015} % Style BST file
\bibliography{references}    % Bibliography file (usually '*.bib' )

\end{document}